\documentclass[a4paper]{article}
\usepackage{amsthm,amsmath,amssymb, bbm}
\usepackage{color}
\newtheorem{teor}{Theorem}[section]
\newtheorem{defin}[teor]{Definition}
\newtheorem{lemm}[teor]{Lemma}
\newtheorem{osse}[teor]{Remark}
\newtheorem{prop}[teor]{Proposition}
\newtheorem{defi}[teor]{Definition}
\newtheorem{coro}[teor]{Corollary}
\newtheorem{prob}[teor]{Problem}
\newtheorem{assu}[teor]{Assumption}

\newcommand{\bele}{\begin{lemm}\begin{sl}}
\newcommand{\enle}{\end{sl}\end{lemm}}
\newcommand{\bedef}{\begin{defi}\begin{sl}}
\newcommand{\eddef}{\end{sl}\end{defi}}
\newcommand{\bete}{\begin{teor}\begin{sl}}
\newcommand{\ente}{\end{sl}\end{teor}}
\newcommand{\beos}{\begin{osse}\begin{rm}}
\newcommand{\eddos}{\end{rm}\end{osse}}
\newcommand{\beas}{\begin{assu}\begin{rm}}
\newcommand{\eddas}{\end{rm}\end{assu}}
\newcommand{\bepr}{\begin{prop}\begin{sl}}
\newcommand{\empr}{\end{sl}\end{prop}}
\newcommand{\bepro}{\begin{prob}\begin{rm}}
\newcommand{\empro}{\end{rm}\end{prob}}
\newcommand{\bede}{\begin{defin}\begin{sl}}
\newcommand{\edde}{\end{sl}\end{defin}}
\newcommand{\beco}{\begin{coro}\begin{sl}}
\newcommand{\enco}{\end{sl}\end{coro}}

\newcommand{\quand}{\quad\text{and}\quad}

\newcommand{\quext}{\quad\text}
\newcommand{\qquext}{\qquad\text}
\newcommand{\de}{\partial}

\newcommand{\RR}{\mathbb{R}}

\newcommand{\NN}{\mathbb{N}}

\newcommand{\beeq}[1]{\begin{equation}\label{#1}}
\newcommand{\eddeq}{\end{equation}}

\newcommand{\beeqa}[1]{\begin{eqnarray}\label{#1}}
\newcommand{\eddeqa}{\end{eqnarray}}

\newcommand{\beal}[1]{\begin{align}\label{#1}}
\newcommand{\eddal}{\end{align}}

\newcommand{\bespl}[1]{\begin{split}\label{#1}}
\newcommand{\edspl}{\end{split}}

\newcommand{\bega}[1]{\begin{gather}\label{#1}}
\newcommand{\edga}{\end{gather}}

\newcommand{\beeqax}{\begin{eqnarray*}}
\newcommand{\eddeqax}{\end{eqnarray*}}

\def\qed{\ifmmode 
  \else \leavevmode\unskip\penalty9999 \hbox{}\nobreak\hfill
  \fi
  \quad\hbox{\hskip.5em\vrule width.4em height.6em depth.05em\hskip.1em}}
\def\endproofsym{\qed}
\renewenvironment{proof}[1][Proof]{\trivlist\item[\hskip\labelsep{\hskip0pt
    {\normalfont\scshape#1.}\hskip .321429\parindent}]\ignorespaces}
{\endproofsym\endtrivlist}
\def\endnobox{\def\endproofsym{}\end{proof}\def\endproofsym{\qed}}

\newcommand{\no}{\nonumber}

\newcommand{\beeqao}{\begin{eqnarray}\no}
\newcommand{\bealo}{\begin{align}\no}
\newcommand{\besplo}{\begin{split}\no}
\newcommand{\begao}{\begin{gather}\no}

\newcommand{\ov}{\overline}

\newcommand{\ba}{\boldsymbol{a}}

\newcommand{\duav}[1]{\langle{#1}\rangle}

\newcommand{\io}{\int_\Omega}

\newcommand{\bV}{\boldsymbol{V}}
\newcommand{\bC}{\boldsymbol{C}}

\newcommand{\epsi}{\varepsilon}
\newcommand{\ee}{_{\varepsilon}}

\newcommand{\OO}{_{\Omega}}
\newcommand{\oo}{_{\Omega}}

\newcommand{\bxi}{\boldsymbol{\xi}}

\newcommand{\bH}{\boldsymbol{H}}

\newcommand{\bn}{\boldsymbol{n}}

\newcommand{\bv}{\boldsymbol{v}}
\newcommand{\dn}{\partial_{\bn}}
\newcommand{\fhi}{\varphi}

\newcommand{\bu}{\boldsymbol{u}}
\newcommand{\vu}{\boldsymbol{u}}

\newcommand{\lhs}{left-hand side}
\newcommand{\rhs}{right-hand side}

\DeclareMathOperator{\dive}{div}
\DeclareMathOperator{\deriv}{d}

\DeclareMathOperator{\Id}{Id}

\DeclareMathOperator{\sign}{sign}

\DeclareMathOperator{\loc}{loc}

\DeclareMathOperator{\Lip}{Lip}

\newcommand{\LDH}{L^2(0,T;H)}
\newcommand{\LDV}{L^2(0,T;V)}
\newcommand{\LDVp}{L^2(0,T;V')}
\newcommand{\LIH}{L^\infty(0,T;H)}
\newcommand{\LIV}{L^\infty(0,T;V)}

\newcommand{\LDHD}{L^2(0,T;H^2(\Omega))}

\let\TeXchi\chi
\def\chi{{\setbox0 \hbox{\mathsurround0pt
$\TeXchi$}\hbox{\raise\dp0 \copy0 }}}

\newcommand{\gammaciapo}{\widehat{\gamma}}

\newcommand{\betaciapo}{\widehat{\beta}}

\newcommand{\calE}{{\mathcal E}}

\newcommand{\dit}{\deriv\!t}

\newcommand{\dir}{\deriv\!r}

\newcommand{\ddt}{\frac{\deriv\!{}}{\dit}}

\newenvironment{giuliorev}{\color{red}}{\color{black}}
\newcommand{\III}{\begin{giuliorev}}
\newcommand{\EEE}{\end{giuliorev}}

\numberwithin{equation}{section}
\begin{document}

\title{On a modified Cahn-Hilliard-Brinkman model with chemotaxis 
and nonlinear sensitivity}

\author{Giulio Schimperna\\
Dipartimento di Matematica, Universit\`a di Pavia\\
and IMATI-CNR, Pavia,\\
Via Ferrata~5, I-27100 Pavia, Italy\\
E-mail: {\tt giulio.schimperna@unipv.it}
}


\maketitle
\begin{abstract}
 We consider an evolutionary PDE system coupling the Cahn-Hilliard equation
 with singular potential, mass source and transport effects,
 to a Brinkman-type relation for the macroscopic velocity field and to a further
 equation describing the evolution of the concentration of a chemical substance
 affecting the phase separation process. The main application we have in mind
 refers to tumor growth models: in particular, the equation for the chemical
 prescribes that such a substance tends to migrate towards the regions where 
 the tumor cells are more dense and consume it more actively.
 The cross-diffusion effects characterizing the system are 
 similar to those occurring in the Keller-Segel model for chemotaxis.
 There is, however, a profound difference between the two settings:
 actually, the Cahn-Hilliard system prescribes a fourth-order dynamics with respect 
 to space variables, whereas most models for chemotaxis are of the 
 second order in space. This fact has a number of specific consequences regarding 
 regularity properties of solutions and conditions ensuring existence.
 Our main results are devoted to proving existence
 of weak solutions in the case when the chemotactic sensitivity function
 depends nonlinearly on the chemical species concentration, and, more
 precisely, has a slow growth at infinity so to avoid finite-time blowup.
 We also analyze the asymptotic problem obtained by letting 
 the viscosity go to zero so to get a Darcy flow regime in the limit.
\end{abstract}

\noindent {\bf Key words:}~~Cahn-Hilliard, chemotaxis, mass source, 
singular potential, nonlinear sensitivity.

\vspace{2mm}

\noindent {\bf AMS (MOS) subject clas\-si\-fi\-ca\-tion:}%
~~35D30, 35K35, 35K86, 35Q92, 92C17, 92C50.

\vspace{2mm}


\section{Introduction}
\label{sec:intro}

Letting $\Omega\subset \RR^d$, $d\in\{2,3\}$, be a smooth bounded domain of 
boundary $\Gamma$, we consider the following PDE system in
$\Omega\times(0,T)$, $T>0$ being an arbitrary, but otherwise fixed, final time
of the evolution process:
\begin{align}\label{CH1}
  & \fhi_t + \vu\cdot \nabla \fhi - \Delta \mu = h(\sigma,\fhi) - \ell \fhi,\\
 \label{CH2}
  & \mu = - \Delta \fhi + f(\fhi) - \chi \sigma,\\
 \label{nutr}
  & \sigma_t + \vu \cdot \nabla \sigma
   - \Delta\sigma + \chi \dive( \alpha(\sigma) \nabla \fhi ) = b (\sigma,\fhi),\\
 \label{brink}
  & - \epsi \dive (D \vu) + \vu 
   = \nabla \pi + \mu \nabla \fhi - \chi \fhi \nabla \sigma,\\
 \label{incompr}
  & \dive \vu = 0.
\end{align} 
Relations \eqref{CH1}-\eqref{CH2}
constitute a suitable form of the Cahn-Hilliard-Oono system
\cite{OO1,OO2} for the order parameter $\fhi$, where the \rhs\ of
\eqref{CH1} plays the role of a mass source decomposed
as the sum of a (dominating) linear part $-\ell\fhi$,
where $\ell>0$ is a given constant, and a nonlinear perturbation 
$h(\sigma,\fhi)$. Convection effects are considered, with 
the macroscopic flow velocity $\vu$ 
satisfying, in an incompressible regime,
the Brinkman-type relation \eqref{brink}, $\pi$
denoting the pressure. In \eqref{brink}, 
the expression $D\vu$ stands for the symmetrized
gradient of $\vu$, with the last two terms on the \rhs\ 
playing the role of Korteweg forces whose specific 
expression guarantees the validity of the energy balance
law. We assume, in principle, $\epsi>0$, 
but (in the spirit, e.g., of 
\cite{KS}), we will also study the behavior as $\epsi\searrow 0$
corresponding to a Darcy regime in the limit.

The nonlinear function $f$ is assumed as the derivative 
(or, more precisely, the subdifferential) of a logarithmic potential $F$ of 
Flory-Huggins type, whose expression reads
\begin{equation}\label{Flog}
  F(r) = (1+r)\log(1+r)+(1-r)\log(1-r)-\frac\lambda2 r^2, \quad r\in[-1,1],
  \quad \lambda \ge 0.
\end{equation}
In particular, as is customary for phase-separation models, the order
parameter is normalized so that the pure configurations are represented
by $\fhi=\pm 1$; correspondingly, the minima of $F$ are attained in 
proximity of these states, and they are deeper for larger $\lambda$.

The main application we have in mind for system \eqref{CH1}-\eqref{incompr}
refers to tumor growth processes; in this setting $\fhi$ may represent the local
proportion of active cancer cells; hence the \rhs\ of \eqref{CH1} drives the
evolution of the total tumor mass. Indeed, assuming suitable boundary
conditions (and, in particular, no-flux boundary conditions
for the {\sl chemical potential}\/ $\mu$, which acts as an auxiliary
variable), the mass evolution law is obtained by 
integrating \eqref{CH1} over $\Omega$. This process is influenced 
by the additional variable 
$\sigma$, which represents the concentration of a chemical substance, 
like a nutrient or a drug. As in related models,  
$\sigma$ is assumed to satisfy the second order 
evolutionary equation \eqref{nutr}, whose specific formulation can be 
considered as the main novelty of the present contribution, so
deserving some words of explanation.

First of all, we point out that tumor growth models based on the Cahn-Hilliard 
``diffuse-interface'' description are becoming increasingly popular 
among the scientific community: for mathematical results and a more 
extensive physical background we may quote, with no claim of exhaustivity, 
the recent papers \cite{CLLW,GL1,GL2,GL3,GLSS,HDvdZO,WLFC}, 
the monograph \cite{CL}, and the references therein. In 
most of the quoted works, actually, a different 
expression for the nutrient equation is assumed, which 
may be written, in the simplest case, as
\begin{equation}\label{nutrlin}
  \sigma_t - \Delta\sigma + \chi \Delta \fhi = b (\sigma,\fhi).
\end{equation}
Here, for the sake of clarity, we have taken constant mobility coefficients and 
omitted the transport effects described by the macroscopic velocity.
Both in \eqref{nutrlin} and in \eqref{nutr}, the coefficient $\chi>0$ 
can be interpreted as a transport parameter which drives the cells in the direction 
of the regions where $\sigma$ takes larger values, basically
providing a chemotactic response of the cells with respect to the nutrient,
an effect which is observed in real world situations.
The main difference between \eqref{nutrlin} and \eqref{nutr}
lies in the expression of the cross-diffusion term, which 
in \eqref{nutrlin} depends linearly on $\fhi$. This, in particular,
implies that $\sigma$ may not satisfy a minimum principle, 
even in the case when the \rhs\ source term $b(\sigma,\fhi)$ 
is properly designed (i.e.\ it satisfies suitable sign and growth conditions). 

Since the nonnegativity of $\sigma$ is an expected property
as $\sigma$ plays the role of a concentration,
it was proposed in \cite{RSS} to replace \eqref{nutrlin} with 
\begin{equation}\label{nutrKS}
  \sigma_t - \Delta\sigma + \chi \dive( \sigma \nabla \fhi ) = b (\sigma,\fhi).
\end{equation}
Actually, such an expression preserves the nonnegativity of $\sigma$ and is 
reminiscent of the Keller-Segel model for chemotaxis,
so providing a new type of coupling between the Cahn-Hilliard equation 
and a Keller-Segel-like relation (we observe that
a different connection between the two models was analyzed in \cite{EPP}).
One of the main notable features of the model introduced
in \cite{RSS} stands in the fact that, differently from the 
``genuine'' Keller-Segel case, the coercivity of 
the energy functional is tied to the uniform boundedness of $\fhi$, hence on
the choice of a singular configuration potential like \eqref{Flog}. 
It is clear, however, that the quadratic behavior of the cross-diffusion term 
in \eqref{nutrKS} may still give rise to blow-up effects; in order to prevent the 
latter, in \cite{RSS} $b$ is assumed to have a logistic growth with 
respect to $\sigma$. This choice, which is rather common in the Keller-Segel
setting, is aimed at penalizing the large values of the concentration;
mathematically speaking, it helps for the sake of obtaining global
in time estimates.

In this paper, considering system \eqref{CH1}-\eqref{incompr},
we will see that, suitably designing the expression 
of the function $\alpha$ in \eqref{nutr},
we may simultaneously keep the minimum principle for $\sigma$
and avoid finite-time blow up under very general and natural 
conditions (and, in particular, with no need for a logistic
growth of $b$). To be precise, we will assume $\alpha$ 
to degenerate at $0$ as fast as $\sigma$, 
and to behave at infinity like a suitable power $\sigma^a$,
where $a\in[0,1)$ is assumed to be small enough in order to
avoid blow up. This choice, usually noted as ``degenerate sensitivity''
in the Keller-Segel literature, is very popular in the chemotaxis community:
one can refer to the seminal paper \cite{HW}
(see also \cite{WWZ} and the references therein for further results 
and a survey of the more recent literature).
Not surprisingly, as happens in similar models characterized by a 
degenerate sensitivity, our main results, which are devoted to 
proving existence and regularity of weak solutions for system 
\eqref{CH1}-\eqref{incompr}, will require suitable restrictions
on the exponent $a$, also depending on the space dimension.

In order to explain the ``spirit'' of our results and to 
describe the novelties and the approach of the present work, 
it is first worth outlining some specific 
features of system \eqref{CH1}-\eqref{nutr} compared
to other models characterized by similar cross-diffusion terms. 
In particular, we may observe the following facts:
\begin{itemize}
 \item The (fourth order) Cahn-Hilliard coupling provides lower time regularity
 but better space regularity for $\fhi$; this effect results in a somehow different 
 behavior of the cross-diffusion term compared to the Keller-Segel model;
 \item The choice of the singular potential \eqref{Flog} guarantees ``for free''
 the uniform boundedness of $\fhi$. This fact not only, as already observed,
 provides coercivity of the energy, but it also guarantees
 additional regularity of the coupling terms; in particular, the uniform 
 boundedness of $\fhi$ is a fundamental information as we seek for 
 additional regularity estimates; 
 \item The Cahn-Hilliard system, which is a subset of \eqref{CH1}-\eqref{incompr}, 
 somehow drives the optimal strategy for getting a-priori estimates.
 Not surprisingly, the argument used for proving existence is primarily based
 on the energy balance, which holds as a direct consequence of the variational structure
 of the system. Then, our assumption $\alpha(\sigma)\sim \sigma^a$
 for large $\sigma$ corresponds to the regularity condition 
 $\sigma\in L^p(\Omega)$, where $p=2-a$; namely,
 a slower growth of $\alpha$ at infinity corresponds to a greater summability
 of solutions. Correspondingly, for large $a$ (or, equivalently,
 small $p$), the outcome of the energy estimate in terms of a-priori regularity
 seems not sufficient to pass to the limit in the approximation, 
 with the main difficulties arising in connection with the product terms and, 
 more specifically, the transport and Korteweg terms depending on $\sigma$;
 \item Looking for additional a-priori regularity, we derive further
 estimates of the so-called ``entropy'' type, 
 according to the terminology currently in use for fourth order evolutionary 
 systems. Basically, this consists in testing \eqref{CH2} by 
 $-\Delta \fhi$ and \eqref{nutr} by $\sigma^{q-1}$ for suitable 
 $q\in[p,\infty)$. Actually, the information resulting from this type of procedure
 provides {\sl separate}\/ regularity bounds for the diffusion
 terms in \eqref{nutr} (and not just for their sum). Moreover, it also 
 yields a  $L^2$-estimate for $\Delta \fhi$ as well
 as some additional summability of $\sigma$; 
 \item Our main result will be devoted to the three-dimensional case:
 we will show that, assuming that $p=2-a\in(12/11,2]$, if the initial datum
 $\sigma_0$ lies in $L^q(\Omega)$ for some $q\in[p,2]$, then
 $\sigma$ keeps staying in $L^q(\Omega)$ uniformly on $(0,T)$.
 In dimension $d=2$, with a simpler proof we will show that 
 similar properties hold for any $p\in(1,2]$. However, the two-dimensional
 case will be analyzed more carefully in a forthcoming paper,
 where we will also discuss the limit case $p=1$. We finally
 observe that the (supposedly) critical exponent $p=12/11$ for 
 $d=3$ seems less restrictive compared to the genuine Keller-Segel setting
 (cf., e.g., \cite{W10});
 \item Considering the presence of (mass and nutrient) source terms,
 we do not expect eventual uniform boundedness of solutions. Indeed, even
 at the level of the energy estimate, our procedure relies in an essential
 way on Gr\"onwall's lemma, implying that solutions may grow exponentially
 in suitable norms. The problem of eventual boundedness of solutions
 (in the spirit of \cite{HW} and several other papers dealing with
 the Keller-Segel model) as well as other questions related with the 
 long-time behavior (e.g., dissipativity, existence of attractors) 
 may be addressed in some future work referring to the case with no 
 external sources, and possibly neglecting the effects of the macroscopic
 velocity;
 \item We finally remark that, in view of the highly nonlinear
 character of the model, uniqueness in the class of weak solutions is 
 not expected to hold. In particular, the combination of a nonlinear
 mass source term in \eqref{CH1} with the singular potential \eqref{Flog},
 as observed in \cite{FS,GGM}, may prevent using a contractive argument;
 in addition to that, the transport effects, especially in the Darcy limit 
 regime, may be difficult to manage.
 Hence, also the question of uniqueness may be considered in a future
 work by assuming a simplified setting.
\end{itemize}

\smallskip 

\noindent%
The plan of the paper is as follows: in the next section, we detail
our assumptions on coefficients and data and state our mathematical 
results. In the subsequent Section~\ref{sec:energy}, we prove
the basic ``physical'' a priori estimates resulting from the 
mass and energy conservation principles and representing the fundamental 
step in the proof of existence. In Section~\ref{sec:apriori2a} we 
obtain further estimates of the so-called ``entropy'' type
which are crucial in order to deal with the cross-diffusion term
in \eqref{nutr} and to provide some further a-priori information
on weak solutions. In Section~\ref{sec:lim} we show weak sequential 
stability of families of weak solutions, namely we prove that 
any sequence of weak solutions complying with the a-priori estimates
uniformly with respect to approximation or regularization
parameters admits at least one limit point which is still
a weak solution. Finally, in Section~\ref{sec:appro} we propose a 
specific regularization scheme and show its compatibility
both with the estimates and with the argument used to pass to the limit.


\section{Assumptions and main results}
\label{sec:main}

We start with introducing a set of notation which will be useful in
order to rigorously formulate our mathematical results. 
Letting $\Omega$ be a smooth bounded domain of $\RR^d$, $d\in\{2,3\}$, 
with boundary $\Gamma$, we 
set $H := L^2(\Omega)$ and $V := H^1(\Omega)$.
We will generally write $H$ in place of $H^d$
(with similar notation for other spaces), whenever vector-valued 
functions are considered. We denote by $(\cdot,\cdot)$ the
standard scalar product of~$H$ and by 
$\| \cdot \|$ the associated Hilbert norm.  Moreover, we equip $V$ 
with the usual norm $\|\cdot\|_V^2 = \|\cdot\|^2 + \|\nabla \cdot\|^2$.
Identifying $H$ with its dual space $H'$ by means of the scalar 
product introduced above, we obtain the chain of continuous 
and dense embeddings $V\subset H \subset V'$.
We will indicate by $\duav{\cdot,\cdot}$
the duality pairing between $V'$ and $V$, or, more generally,
between $X'$ and $X$, where $X$ is a generic Banach space continuously
and densely embedded into $H$.

Next, we assume $F$ be given by \eqref{Flog}, we indicate as $f$ the 
derivative (or, more precisely, the subdifferential) of $F$, and we 
denote as $\beta(r):=f(r)+\lambda r$ the ``monotone part''
of $f$, namely
\begin{equation}\label{defi:beta}
  \beta(r) = \log(1+r) - \log(1-r), \quad 
   r \in (-1,1).
\end{equation}
We also set 
\begin{equation}\label{defi:bV}
  \bV:=\big\{\ba \in V^d,\,~\dive \ba \equiv 0\,~\text{in }\,\Omega,%
  \,~\ba\cdot\bn = 0\,~\text{on }\,\Gamma\big\},
\end{equation}
with $\bn$ denoting the outer normal unit vector to $\Gamma$.
We also let $H^2_{\bn}(\Omega)$ be the (closed) subspace of $H^2(\Omega)$ 
containing the functions with zero normal derivative on $\Gamma$.
For $\bv\in\bV$, we denote as $D\bv=\frac{\nabla \bv + (\nabla\bv)^t}2$ 
the symmetrized gradient of $\bv$.
As in \cite{GLRS} (see also \cite{FLRS} for a ``multi-phase'' version of 
this condition), we assume that 
\begin{equation}\label{hp:h}
  h\in C^1(\RR\times\RR), \qquad 
   | h(\sigma,\fhi) | \le H \quext{for every }\,(\sigma,\fhi)\in [0,\infty)\times\RR,
\end{equation}
for a suitable constant $H>0$ satisfying the compatibility condition 
\begin{equation}\label{comp:h}
   \frac{H}{\ell} < 1.
\end{equation}
Note that the above implies in particular that, uniformly with respect to $\sigma$, 
$h(\sigma,\fhi) - \ell \fhi$ is strictly negative (respectively, strictly positive) 
in a left neighbourhood of $\fhi=1$ (respectively, in a right neighbourhood 
of $\fhi=-1$). In addition, we ask that there exist $C_{h,1},C_{h,2}>0$ such that the 
partial derivatives of $h$ satisfy
\begin{equation}\label{hp:h2}
   | h_\fhi(\sigma,\fhi) | \le C_{h,1}, \quand  
   | h_\sigma(\sigma,\fhi) | \le C_{h,2} \frac1{1+|\sigma|},
   \quext{for every }\,(\sigma,\fhi) \in[0,\infty)\times\RR.   
\end{equation}
The second condition, which has a mainly technical character, corresponds 
to asking that, if the concentration $\sigma$ is very large, 
its effects on the total mass tend to become essentially indepdendent of 
variations of $\sigma$. Actually, this assumption, which
seems reasonable from the physical viewpoint, is probably not optimal
mathematically: the proof could be likely
adaptable to the case when $h_\sigma$ behaves at 
infinity as $\sigma^{-k}$ for some suitable $k\in(0,1)$ depending on
$p$; however, since this does not appear to be an essential
point, we decided to assume a simpler condition so 
to reduce technicalities. We also observe that, in \eqref{hp:h},
\eqref{hp:h2} (and elsewhere below) we considered the possibility
that $|\fhi|>1$: actually, while ``in the limit'' we will 
have $|\fhi|\le 1$ almost everywhere due to the occurrence of 
the singular potential, this may not be the case in an approximation
as $F$ is smoothed out (see Section~\ref{sec:appro} below for details).
On the other hand, the assumptions are only stated 
for $\sigma \in [0,\infty)$ as it is clear that also in the approximation
the minimum principle for $\sigma$ is satisfied.

As noted in the introduction, we assume that
$\alpha(\sigma)\sim\sigma$ near zero and 
$\alpha(\sigma)\sim\sigma^a$, $a\in[0,1)$ near infinity (here and below,
if $a=0$, of course, we interpret $\sigma^a \equiv 1$ for $\sigma\ge 0$).
We also set $p = 2-a$, so that $p\in(1,2]$,
in order to emphasize the $L^p$-summability provided by the above ansatz.
Then, using $p$, our assumption on $\alpha$ can be stated as
\begin{equation}\label{hp:alpha}
  \alpha(s) = \frac{s}{1+s^{p-1}}, 
   \quext{for }\,s\ge 0.
\end{equation}
Consequently, we have
\begin{equation}\label{alpha2}
  \frac{1}{\alpha(s)} = \frac{1+s^{p-1}}s
   = \frac1s + s^{p-2},
   \quext{for }\,s>0.
\end{equation}
We can then define
\begin{equation}\label{gamma}
   \gamma(s) := \int_1^s \frac{\dir}{\alpha(r)}
    = \ln s + \frac{1}{p-1} s^{p-1} - \frac{1}{p-1},
    \quext{for }\,s>0.
\end{equation}
%
%
With the above notation at hand, equation \eqref{nutr} can be more 
conveniently restated as
\begin{equation}\label{nutr2}
  \sigma_t + \vu \cdot \nabla \sigma
   - \dive \big( \alpha(\sigma) \nabla ( \gamma(\sigma) - \chi \fhi ) \big) = b (\sigma,\fhi).
\end{equation}
%
%
%
Next, we assume the nutrient source term to satisfy
%
%
\begin{equation}\label{hp:b}
   b\in \Lip_{\loc}([0,\infty)\times\RR;\RR), \qquad
    -b_0 \sigma \le b(\sigma,\fhi) \le b_\infty ( 1 + \sigma),
    \quext{ for every }\,(\sigma,\fhi)\in  [0,\infty)\times \RR,
\end{equation}
and for some constants $b_0,b_\infty>0$. Note that the above condition
is designed so to ensure the applicability of a minimum principle argument
for $\sigma$ uniformly with respect to $\fhi$. 

System \eqref{CH1}-\eqref{incompr} is complemented with the Cauchy conditions 
\begin{equation}\label{init}
  \fhi|_{t=0} = \fhi_0, \qquad
  \sigma|_{t=0} = \sigma_0,
   \quext{in }\,\Omega,
\end{equation}
with the initial data satisfying at least the following properties:
\begin{align}\label{hp:fhi0}
  & \fhi_0 \in V, \quad F(\fhi_0)\in L^1(\Omega),
   \quad m_0:=(\fhi_0)\OO \in (-1,1),\\
 \label{hp:sigma0}
   & \sigma_0 \in L^p(\Omega), \quad \sigma_0>0~~\text{a.e.\ in }\, \Omega,
   \quad \ln \sigma_0 \in L^1(\Omega).
\end{align}
The above conditions correspond to the finiteness of the initial energy
and to the fact that, as is customary in the setting of the Cahn-Hilliard
system with logarithmic potential \eqref{Flog}, we cannot admit inital configurations 
where $\fhi$ is almost everywhere equal to $1$ (or to $-1$). 
Here and below, for a generic function $v$ defined of $\Omega$
we have denoted by $v\OO$ the spatial mean of $v$,
namely
\begin{equation}\label{mean}
  v\OO := \frac{1}{|\Omega|} \io v,
\end{equation}
$|\Omega|$ denoting the $d$-dimensional Lebesgue measure of $\Omega$.

Moreover we shall assume no-flux (i.e., homogeneous Neumann) boundary conditions 
for $\mu$, $\fhi$ and $\sigma$, and {\sl complete slip}\/
boundary conditions for the velocity,
namely
\begin{align}\label{noflux}
  & \dn \mu = \dn \fhi = \dn \sigma = 0,
   \qquext{on }\,\Gamma\times(0,T),\\
 \label{navier}
  & \bu\cdot \bn = 0, \quad
  - \epsi( ( D\vu ) \bn)_\tau 
       = 0, 
   \qquext{on }\,\Gamma\times(0,T).
\end{align}
Here and below, for a (sufficiently smooth) vector field $\ba$ defined on $\Gamma$, 
we have denoted by $\ba_\tau:=\ba - (\ba \cdot \bn)\bn$ the tangential component
of $\ba$.

\smallskip

\noindent%
With the above material at hand, we can now present 
our main results, which are devoted to proving 
existence of global in time weak solutions to the 
Cahn-Hilliard-Brinkman-nutrient model as well as to studying 
the asymptotic behavior of solutions as $\epsi\searrow 0$.
As observed in the introduction, we will consider separately the cases 
$d=2$ and $d=3$. In the two-dimensional case we will prove the following 
simpler (but likely not optimal) statement:
\bete\label{teo:main}
 Let the assumptions\/  \eqref{Flog}, \eqref{hp:h}-\eqref{hp:h2}, \eqref{hp:alpha},
 \eqref{hp:b}, \eqref{hp:fhi0}-\eqref{hp:sigma0} hold and let $d=2$. Let $\epsi>0$,
 and assume that $p>1$.
 %
 %
 Moreover, let also
 \begin{equation}\label{sigma0:en}
   \sigma_0 \in H.
 \end{equation}
 Then, there exists at least one quadruple $(\fhi,\mu,\sigma,\vu)$ with the regularity
 properties
 \begin{align}\label{rego:fhi}
  & \fhi \in H^1(0,T;V')\cap \LIV \cap L^4(0,T;H^2(\Omega))
    \cap L^2(0,T;W^{2,P}(\Omega)),\\
  \label{rego:Ffhi}
  & F(\fhi) \in L^\infty(0,T;L^1(\Omega)),
   \qquad f(\fhi) \in L^2(0,T;L^{P}(\Omega)),\\
 \label{rego:mu}
  & \mu\in \LDV,\\
 \label{rego:vu}
  & \vu \in L^2(0,T;\bV),\\
 \label{rego:sigma}
  & \sigma \in H^1(0,T;W^{1,R}(\Omega)') 
    \cap L^\infty(0,T;H) \cap L^2(0,T;V),\\
 \label{pos:sigma}
  & \sigma>0~~\text{a.e.\ in }\,\Omega\times(0,T), \qquad
   \ln \sigma \in \LDV,\\
 \label{rego:H}
  & \bH(\sigma,\fhi) := \alpha^{1/2}(\sigma)\nabla(\gamma(\sigma)-\chi\fhi) \in \LDH,
 \end{align}
 where the exponents $P$ in\/ \eqref{rego:fhi}-\eqref{rego:Ffhi} and $R$ in\/
 \eqref{rego:sigma} may be arbitrarily taken in $[2,\infty)$ and in 
 $(2,\infty]$, respectively. 
 The quadruple $(\fhi,\mu,\sigma,\vu)$ satisfies, a.e.~in $(0,T)$, the following 
 weak version of system \eqref{CH1}-\eqref{incompr}:
 \begin{align}\label{CH1w}
   & \duav{\fhi_t,\xi} - \io \fhi\vu\cdot\nabla \xi
    + \io \nabla \mu\cdot\nabla \xi = \io \big( h(\sigma,\fhi) - \ell \fhi\big)\xi
     \quext{for every }\,\xi\in V,\\
  \label{CH2w}
   & \mu = - \Delta \fhi + f(\fhi) - \chi \sigma,
   \quext{a.e.\ in }\,\Omega,\\
  \label{nutrw}
   & \duav{\sigma_t,\eta} - \io \sigma\vu \cdot \nabla \eta
    + \io \alpha^{1/2}(\sigma) \bH(\sigma,\fhi) \cdot \nabla \eta 
    = \io b (\sigma,\fhi) \eta,
    \quext{for all }\,\eta\in W^{1,R}(\Omega),\\
  \label{brinkw}
   & \epsi \io (D \vu) : ( D \bxi) + \io \vu \cdot \bxi
    = - \io \fhi\bxi \cdot \nabla \mu 
    + \chi \io \sigma \bxi\cdot \nabla\fhi
    \quext{for every }\,\bxi\in \bV.
 \end{align} 
 Moreover, there hold the initial conditions \eqref{init}
 as well as the boundary condition
 \begin{equation}\label{bound}
   \dn \fhi = 0
   \quext{a.e.\ on }\,\Gamma\times(0,T).
 \end{equation}
 Finally, let, for $\epsi\in(0,1)$, $(\fhi\ee,\mu\ee,\sigma\ee,\vu\ee)$ be 
 a family of solutions satisfying the system in the sense specified above.
 Then, there exist a (nonrelabelled) subsequence of $\epsi\searrow0$ and 
 a limit quadruple $(\fhi,\mu,\sigma,\vu)$ such that 
 \begin{align}\label{coee:11}
   & \fhi\ee \to \fhi \quext{weakly star in }\,
    H^1(0,T;V')\cap \LIV \cap L^4(0,T;H^2(\Omega)),\\
  \label{coee:12}
   & \mu\ee \to \mu \quext{weakly in }\,\LDV,\\
  \label{coee:13}
   & \sigma\ee \to \sigma \quext{weakly star in }\,
    W^{1,4/3}(0,T;W^{1,4}(\Omega)') \cap L^\infty(0,T;H) \cap L^2(0,T;V),\\
  \label{coee:14}
   & \vu\ee \to \vu \quext{weakly in }\,\LDH.
 \end{align} 
 The quadruple $(\fhi,\mu,\sigma,\vu)$ satisfies\/ \eqref{CH1w}-\eqref{brinkw}
 with $\epsi=0$ (where, now, $R=4$, i.e., $\eta$ has to lie in $W^{1,4}(\Omega)$,
 in \eqref{nutrw}), the additional boundary condition \eqref{bound}, and the 
 initial conditions \eqref{init}. 
\ente
\beos\label{rem:upress}
 In the Darcy limit $\epsi=0$, the velocity $\vu$ just acts as an auxiliary 
 variable and the system could be restated directly in terms of the pressure 
 $\pi$ (or, more precisely, of its gradient as $\pi$ is determined
 up to an additive constant). Actually, $\pi$ turns out to satisfy
 the limit elliptic system
 \begin{equation}\label{brilim}
   - \Delta \pi = 
    \dive( \mu \nabla \fhi - \chi \fhi \nabla \sigma),
 \end{equation}
 complemented with the natural (no-flux) boundary conditions resulting from
 a comparison with \eqref{noflux}. Once $\pi$ is given by \eqref{brilim}, 
 \eqref{brink} with $\epsi=0$ can then
 be seen as a {\sl definition}\/ of $\vu$ in terms of the other variables.
\eddos
\noindent%
We now move to our main statement referring to the three-dimensional
setting. In this case, we will first prove that, if $p$ is strictly larger than 
the (supposedly) critical exponent $12/11$, then the two summands in the 
cross-diffusion term $\bH(\sigma,\fhi)$ can be controlled {\sl separately};
this information, which is achieved by means of the so-called ``entropy'' method 
for fourth order PDE's, is crucial for the sake of proving existence 
of weak solutions. Furthermore, using a more refined version of the entropy
estimate, we can also show that, if $\sigma_0\in L^q(\Omega)$ for $q\in[p,2]$, 
then a natural parabolic estimate for $\sigma^{q/2}$ holds (see \eqref{rego:sigma1b}
below). Finally, if $q>6/5$, we can take the limit $\epsi\to 0$
and prove convergence to the Darcy flow. This further restriction 
is motivated by the need for having sufficient regularity 
to control, uniformly as $\epsi\searrow0$, the transport and Korteweg terms
depending on $\sigma$. These results are summarized in the following statement:
\bete\label{teo:main2}
 Let the assumptions\/  \eqref{Flog}, \eqref{hp:h}-\eqref{hp:h2}, \eqref{hp:alpha},
 \eqref{hp:b}, \eqref{hp:fhi0}-\eqref{hp:sigma0} hold, and let $d=3$.
 Let, moreover, $\epsi>0$ and let
 \begin{equation}\label{hp:p2}
    p\in(12/11,2].
 \end{equation}
 In addition to that, let us additionally assume 
 \begin{equation}\label{sigma0:en2}
   \sigma_0 \in L^q(\Omega) \quext{where }\,q\in [p,2].
 \end{equation}
 Then, there exists at least one quadruple $(\fhi,\mu,\sigma,\vu)$ with the regularity
 properties
 \begin{align}\label{rego:fhi1}
  & \fhi \in H^1(0,T;V')\cap \LIV 
    \cap L^{P_0}(0,T;H^2(\Omega)) \cap L^2(0,T,W^{2,\frac{3q}{3-q}}(\Omega)),\\
 \label{rego:Ffhi1}
  & F(\fhi) \in L^\infty(0,T;L^1(\Omega)),
   \qquad f(\fhi) \in L^2(0,T;L^{\frac{3q}{3-q}}(\Omega)), \\
 \label{rego:mu1}
  & \mu\in \LDV,\\ 
 \label{rego:vu1}
  & \vu \in L^2(0,T;\bV),\\
 \label{rego:sigma1a}
  & \sigma \in W^{1,Z}(0,T;W^{1,R}(\Omega)') \cap L^2(0,T;W^{1,S}(\Omega))
   \quext{for some }Z>1,\\
 \label{rego:sigma1b}
  & \sigma^{q/2} \in L^\infty(0,T;H) \cap L^2(0,T;V),\\
 \label{pos:sigma1}
  & \sigma>0~~\text{a.e.\ in }\,\Omega\times(0,T), \qquad
   \ln \sigma \in \LDV,\\
 \label{rego:H1}
  & \bH(\sigma,\fhi) := \alpha^{1/2}(\sigma)\nabla(\gamma(\sigma)-\chi\fhi) \in \LDH,
 \end{align}
 where the exponents $P_0,R,S$ in \eqref{rego:fhi1}, \eqref{rego:sigma1a} 
 satisfy (note that $S>1$ due to \eqref{hp:p2})
 \begin{equation}\label{P0}
   P_0 = \min \Big\{ \frac{18q-6p}{12-5p}, 4 \Big\}, 
    \qquad
   S = \min \Big\{ \frac{6p}{12-5p}, p \Big\},
   \qquad
   R = \max \big\{ 4, p' \big\},
 \end{equation}
 $p'$ being the conjugate exponent to $p$. The quadruple $(\fhi,\mu,\sigma,\vu)$
 satisfies, a.e.~in $(0,T)$, relations \eqref{CH1w}-\eqref{brinkw}, where, now,
 the exponent $R$ in \eqref{nutrw} is specified by \eqref{P0}.
 Moreover, there hold the boundary condition~\eqref{bound}
 and, in the sense of traces, the initial conditions~\eqref{init}.\\
 Next, let us additionally assume 
 \begin{equation}\label{su:q}
    q > 6/5
 \end{equation}
 and let, for $\epsi\in(0,1)$, $(\fhi\ee,\mu\ee,\sigma\ee,\vu\ee)$ be 
 a family of solutions satisfying the system in the sense specified above.
 Then, there exist a (nonrelabelled) subsequence of $\epsi\searrow0$ and 
 a limit quadruple $(\fhi,\mu,\sigma,\vu)$ such that 
 \begin{align}\label{coee:21a}
   & \fhi\ee \to \fhi \quext{weakly in }\,
    H^1(0,T;V')\cap L^{P_0}(0,T;H^2(\Omega)) \cap L^2(0,T,W^{2,\frac{3q}{3-q}}(\Omega)),\\
   \label{coee:21b}
   & \fhi\ee \to \fhi \quext{weakly star in }\, \LIV,\\
  \label{coee:22}
   & \mu\ee \to \mu \quext{weakly in }\,\LDV,\\
  \label{coee:23}
   & \sigma\ee \to \sigma \quext{weakly in }\,
    W^{1,Z}(0,T;W^{1,Z'}(\Omega)') \cap L^2(0,T;W^{1,S}(\Omega)), 
    \quext{for some }\,Z>1,\\
  \label{coee:24}
   & \vu\ee \to \vu \quext{weakly in }\,\LDH.
 \end{align} 
 The quadruple $(\fhi,\mu,\sigma,\vu)$ satisfies\/ \eqref{CH1w}-\eqref{brinkw}
 with $\epsi=0$ (where, now, in \eqref{nutrw}, $\eta$ has to lie in $W^{1,Z'}(\Omega)$,
 $Z'$ being the ``large'' conjugate exponent to $Z$), the additional boundary
 condition \eqref{bound}, and the initial conditions \eqref{init}. 
\ente
\noindent%
We conclude this section with some remarks aimed at clarifying some aspects 
of the above statements.
\beos\label{rem:teomain}
It is worth observing that the boundary conditions for $\mu$ and $\vu$
are embedded, respectively, in the weak formulation \eqref{CH1w}
and in the regularity \eqref{rego:vu} (or \eqref{rego:vu1})
combined with relation \eqref{brinkw}.
We also remark that relation \eqref{nutrw} might be more standardly 
rewritten as 
\begin{equation}\label{nutrw2}
  \duav{\sigma_t,\eta} - \io \sigma\vu \cdot \nabla \eta
    + \io \alpha(\sigma) \nabla(\gamma(\sigma)-\chi\fhi) \cdot \nabla \eta 
    = \io b (\sigma,\fhi) \eta.
\end{equation}
We preferred to keep formulation \eqref{nutrw} since it better emphasizes
the additional regularity \eqref{rego:H} satisfied by the mixed term $\bH$.
In particular we point out that, while $\bH$, as a sum, lies in $L^2$,
as it is decomposed into the summands $\alpha^{1/2}(\sigma) \nabla \gamma(\sigma)$
and $-\chi \alpha^{1/2}(\sigma) \nabla \fhi$, these summands may
fulfill lower regularity properties.
\eddos
\beos\label{rem:teomain0}
Referring to the statement of Theorem~\ref{teo:main2}, some comments about assumptions 
\eqref{hp:p2} and \eqref{sigma0:en2} are in order. First of all, the condition $p>12/11$ 
corresponds to the most general situation in which we are able to prove, via the entropy
method, decoupled regularity for the two summands in the cross-diffusion
term $\bH(\sigma,\fhi)$. In particular, one may well consider the case $q=p$,
which corresponds to proving existence for initial data having exactly the 
``energy'' regularity. On the other hand, looking at \eqref{sigma0:en2},
one may assume that $\sigma_0$ enjoys some additional summability; in that case, 
we can prove that this level of summability is maintained in time 
(cf.\ \eqref{rego:sigma1b}).
\eddos
\beos\label{rem:teomain2}
We decided not to treat the case when $\sigma_0 \in L^q(\Omega)$ for $q$ 
{\sl strictly}\/ larger than $2$. In that situation,
however, it is expected that the regularity
of solutions could be further improved by bootstrapping. Moreover, assuming some
condition on $\nabla\sigma_0$ one may also prove some summability of
$\sigma_t$ and $\Delta\sigma$ (which are now  controlled only in
Sobolev spaces of negative order with respect to space
variables). In turn, this may help improving
the information on $\fhi$ by means of the so-called ``second energy estimate''
that can be performed for the Cahn-Hilliard equation with singular potential.
The question of additional regularity of solutions may be treated 
in a future work, possibly neglecting, also in this case,
the effects of convection.
\eddos
%
%
%
%

%


\section{Mass balance and energy estimate}
\label{sec:energy}

In this part we derive a number of a priori estimates satisfied by any hypothetical
solution of our system and holding under the natural ``physical'' conditions
on the initial data, namely finiteness of the physical energy, compatibility of the
initial mass with the configuration potential,
and positivity of the initial concentration. For simplicity, we will work
directly on the ``original'' system \eqref{CH1}-\eqref{incompr}: in Section~\ref{sec:appro}
we will see how the procedure may be adapted to a specific approximation scheme.
We start with deriving the evolution of mass.

\smallskip

\noindent%
{\bf Balance of mass.}~~%
Integrating \eqref{CH1} over $\Omega$ and using the incompressibility constraint
and the boundary conditions, we immediately deduce 
\begin{equation}\label{mass}
  \ddt \fhi\OO + \ell\fhi\OO =
    \frac1{|\Omega|} \io h(\sigma,\fhi).
\end{equation}
Applying \eqref{comp:h} and using the last \eqref{hp:fhi0} with a simple argument
based on the comparison principle for ODE's, we deduce that there exists a 
sufficiently small $\delta\in (0,1)$, depending only on $m_0$, $\ell$ and $H$ (and in 
particular independent of $T$) and such that 
\begin{equation}\label{medie}
  -1 + \delta \le m(t) = (\fhi(t))\OO \le 1-\delta 
   \quext{for every }\,t\in[0,T].
\end{equation}
%

\noindent%
{\bf Positivity of the nutrient concentration.}~~%
As in similar models,
this property may be obtained by means of a standard argument based on the 
Stampacchia truncation method combined with a regularization of equation
\eqref{nutr}, which is needed in order to justify the use of the truncated
$\sigma$ as a test function. In particular, one may refer to the regularization
scheme that will be detailed in Section~\ref{sec:appro} and possibly operate
a further smoothing of $\vu$ and $\fhi$. Note also that, by \eqref{hp:b}, the 
\rhs\ of \eqref{nutr} is suitably designed so to allow for such an argument.
It is also worth observing that
the spirit of the Stampacchia truncation argument is to prove
\begin{equation}\label{stamp}
   \fhi_0(\cdot)\ge 0~~\text{a.e.\ in }\,\Omega~~\Longrightarrow~~%
    \fhi(\cdot,\cdot)\ge 0~~\text{a.e.\ in }\,\Omega\times(0,T).
\end{equation}
However, we point out that \eqref{stamp} is not sufficient for our purposes; 
actually, in the sequel, based on \eqref{hp:sigma0}, we will improve the above
argument obtaining a stronger information.

\smallskip

\noindent%
{\bf Energy estimate.}~~%
We now prove the most important
a-priori bound, which is a direct consequence of the variational structure 
of the system. To derive it, we test \eqref{CH1} by $\mu$, 
\eqref{CH2} by $\fhi_t$, and take the difference so to obtain
\begin{equation}\label{en:11}
  \ddt \Big(\frac12 \| \nabla \fhi \|^2 + \io F(\fhi) \Big) 
   - \chi \io \sigma \fhi_t
   + \| \nabla\mu \|^2 
   = - \io \mu \vu \cdot \nabla\fhi
   + \io \big( h(\sigma,\fhi) - \ell\fhi \big) \mu.
\end{equation}
Moreover, we test \eqref{nutr2} by $\gamma(\sigma) - \chi\fhi$. This gives 
\begin{equation}\label{en:12}
  \ddt \io \gammaciapo(\sigma)
   - \chi \io \sigma_t \fhi
   + \io \alpha(\sigma) \big| \nabla ( \gamma(\sigma) - \chi \fhi ) \big|^2
  = \io b(\sigma,\fhi) (\gamma(\sigma) - \chi \fhi) 
   + \chi \io \fhi \vu \cdot \nabla \sigma,
\end{equation}
where the incompressibility constraint \eqref{incompr}
and the first \eqref{navier} have also been used.
Note also that $\gamma(\sigma)$ contains $\ln\sigma$ a summand; as will
be clear in the sequel of the argument, $\ln\sigma$ can be shown to 
have some a-priori regularity, which justifies its use as a test
function. Here, we have also set
\begin{equation}\label{gammaciapo}
  \gammaciapo(s) := \int_1^s \gamma(r)\,\dir 
   = \frac1{p(p-1)} s^p + s \ln s - \frac{p}{p-1} s + \frac{p+1}p.
\end{equation}
Testing now \eqref{brink} by $\vu$, and using once more \eqref{incompr}
and \eqref{navier}, we deduce
\begin{equation}\label{en:13}
  \epsi \| D \vu \|^2 
   + \| \vu \|^2
  = \io \mu \vu \cdot \nabla\fhi
   - \chi \io \fhi \vu \cdot \nabla\sigma.
\end{equation}
Summing \eqref{en:11}, \eqref{en:12} and \eqref{en:13}, we obtain
\begin{align}\no 
  & \ddt \Big(\frac12 \| \nabla \fhi \|^2 + \io F(\fhi) 
   + \io \big( \gammaciapo(\sigma) - \chi \sigma\fhi \big) \Big)
   + \| \nabla\mu \|^2 
   + \io \alpha(\sigma) \big| \nabla ( \gamma(\sigma) - \chi \fhi ) \big|^2\\
 \label{en:14}
  & \mbox{}~~~~~
  + \epsi \| D \vu \|^2 
   + \| \vu \|^2
    = \io b(\sigma,\fhi) (\gamma(\sigma) - \chi \fhi)
   + \io \big( h(\sigma,\fhi) - \ell\fhi \big) \mu.
\end{align}
Then, defining now the {\sl energy functional}
\begin{equation}\label{energy}
  \calE(\fhi,\sigma)
   := \frac12 \| \nabla \fhi \|^2 + \io F(\fhi) 
   + \io \big( \gammaciapo(\sigma) - \chi \sigma\fhi \big),
\end{equation}
and recalling \eqref{gammaciapo}, it is immediate to check the coercivity property
\begin{equation}\label{coer:energy}
  \calE(\fhi,\sigma)
   \ge \kappa \big( \| \fhi \|_V^2 + \| F(\fhi) \|_{L^1(\Omega)} \big)
    + \kappa_p \| \sigma \|_{L^p(\Omega)}^p - c_p,
\end{equation}
for suitable constants $\kappa,\kappa_p>0$ and $c_p\ge 0$.
Analogously, one may also prove a control from below, namely
for some $C>0$ there holds
\begin{equation}\label{bdd:energy}
  \calE(\fhi,\sigma)
   \le C \big( \| \fhi \|_V^2 + \| F(\fhi) \|_{L^1(\Omega)} 
        + \| \sigma \|_{L^p(\Omega)}^p + 1 \big).
\end{equation}
Of course, properties \eqref{coer:energy} and \eqref{bdd:energy} are to 
be intended to hold for $(\fhi,\sigma)$ belonging to the natural domain
of the energy functional (namely, for the pairs $(\fhi,\sigma)$ such that
$\fhi\in V$, $F(\fhi)\in L^1(\Omega)$, $\sigma\in L^p(\Omega)$ with 
$\sigma\ge 0$ a.e.\ in $\Omega$).
\beos\label{oss:limi1}
 In order to deduce \eqref{coer:energy} and 
 \eqref{bdd:energy}, we have used that
 \begin{equation}\label{duapp}
   \chi \bigg| \io \sigma\fhi \bigg|
    \le c \| \sigma \|_{L^p(\Omega)} \| \fhi \|_{L^{p'}(\Omega)}
    \le c \| \sigma \|_{L^p(\Omega)},
 \end{equation}
 where $p'$ is the conjugate exponent to $p$ and the last inequality follows
 from the uniform boundedness of $\fhi$, which 
 is a consequence of the choice of the potential \eqref{Flog}.
 It is worth however observing that, in order for the above to work, it would be 
 sufficient to assume that $F$ satisfies the weaker coercivity property
 \begin{equation}\label{co:q0}
   F(\fhi) \ge \kappa \| \fhi \|^{q_0}_{L^{q_0}(\Omega)} - c,
 \end{equation}
 for some $q_0>p'$, where $c,\kappa>0$. Actually, in Section~\ref{sec:appro} 
 below, we will exhibit an approximation of $F$ such that
 \eqref{co:q0} holds with $c,\kappa>0$ independent of the approximation
 parameter (see \eqref{betan11} in the statement of Lemma~\ref{lem:linfty}
 for details, see also \cite{RSS} for additional considerations).
\eddos
\noindent%
In order to control the terms on the \rhs\ of \eqref{en:14}, 
we first observe that, by \eqref{hp:b} and \eqref{gamma},
\begin{align}\no
  \io b(\sigma,\fhi) (\gamma(\sigma) - \chi \fhi)
   & = \int_{\{\sigma < 1\}} b(\sigma,\fhi) (\gamma(\sigma) - \chi \fhi)
    + \int_{\{\sigma\ge 1\}} b(\sigma,\fhi) (\gamma(\sigma) - \chi \fhi)\\
 \no
  & \le c \int_{\{\sigma < 1\}} ( 1 + | \sigma\ln\sigma | )
    + c \int_{\{\sigma\ge 1\}} ( 1 + \sigma^p )\\
 \label{st:b} 
  & \le c \big( 1 + \| \sigma \|_{L^p(\Omega)}^p \big).
\end{align}
Here, we used once more the uniform boundedness of $\fhi$ 
(but there hold the considerations made in Remark~\ref{oss:limi1}).
Next, regarding the term depending on the chemical potential, we 
replace the expression \eqref{CH2} of $\mu$ therein, so to obtain
\begin{align}\no
  \io \big( h(\sigma,\fhi) - \ell\fhi \big) \mu
  & = \io \big( h(\sigma,\fhi) - \ell\fhi \big) 
    \big( - \Delta \fhi + f(\fhi) - \chi\sigma)\\
 \no
  & = - \io \big( h(\sigma,\fhi) - \ell\fhi \big) \Delta \fhi 
   + \io \big( h(\sigma,\fhi) - \ell\fhi \big) f(\fhi) 
   - \io \chi \big( h(\sigma,\fhi) - \ell\fhi \big) \sigma\\
 \label{en:11b}   
  & =: - I_1 + I_2 - I_3.
\end{align}
Let us manage the quantities on the \rhs. First,
let us recall that, as noted above, $h(\sigma,\fhi)-\ell\fhi$ 
has the opposite sign to $\fhi$
(hence to $f(\fhi)$) when $|\fhi|\sim 1$ and uniformly in $\sigma$.
This implies that $I_2\le c$. Next, by 
\eqref{hp:h} with the boundedness of $\fhi$ (recall, however,
Remark~\ref{oss:limi1}), we deduce
\begin{equation}\label{st:h1}
  | I_3 | \le c \| \sigma \|_{L^1(\Omega)}
   \le c \big( 1 + \| \sigma \|_{L^p(\Omega)}^p \big).
\end{equation}
In order to control $I_1$ we observe that, recalling assumption \eqref{hp:h2}
and applying Young's inequality,
\begin{align}\no
  - I_1 & = \io h_\sigma(\sigma,\fhi) \nabla\sigma \cdot \nabla\fhi 
   + \io ( h_{\fhi}(\sigma,\fhi) - \ell ) |\nabla \fhi|^2\\
 \label{st:h2}
  & \le C_{h,2} \io \frac{1}{1+\sigma} |\nabla\sigma| |\nabla\fhi| 
   + c \| \nabla\fhi \|^2
  \le \frac14 \| \nabla\ln \sigma \|^2
   + c \| \nabla\fhi \|^2.
\end{align}
Then, collecting the above considerations, \eqref{en:14} yields  
\begin{align}\no 
  & \ddt \Big(\frac12 \| \nabla \fhi \|^2 + \io F(\fhi) 
   + \io \big( \gammaciapo(\sigma) - \chi \sigma\fhi \big) \Big)
   + \| \nabla\mu \|^2 
   + \io \alpha(\sigma) \big| \nabla ( \gamma(\sigma) - \chi \fhi ) \big|^2\\
 \label{en:14b}
  & \mbox{}~~~~~
  + \epsi \| D \vu \|^2 
   + \| \vu \|^2
   \le c \big( 1 + \| \sigma \|_{L^p(\Omega)}^p \big)
   + c \| \nabla\fhi \|^2
   + \frac14 \| \nabla\ln \sigma \|^2.
\end{align}
The simplest way to control the last term consists in improving a bit the
positivity estimate on $\sigma$: to this aim, we test \eqref{nutr2} 
by $-1/\sigma$. Also this procedure could be justified 
by truncation. We then deduce
\begin{equation}\label{co:24}
  \ddt \io -\ln \sigma 
   + \io \frac{\gamma'(\sigma)\alpha(\sigma)}{\sigma^2} | \nabla \sigma |^2
  = - \io \frac{b(\sigma,\fhi)}{\sigma}
   + \chi \io \frac{\alpha(\sigma)}{\sigma^2} \nabla\fhi \cdot \nabla \sigma.
\end{equation}
Now, owing to \eqref{hp:b} and considering that $\sigma\ge 0$ by the Stampacchia
truncation argument, we deduce that $\sigma^{-1} b(\sigma,\fhi) \le b_0$
almost everywhere. Next, by \eqref{gamma} one has (as expected)
\begin{equation}\label{co:32ba}
  \io \frac{\gamma'(\sigma)\alpha(\sigma)}{\sigma^2} | \nabla \sigma |^2
   = \io \frac{1 + \sigma^{p-1}}{\sigma} \frac{\sigma}{1 + \sigma^{p-1}}
        \frac{1}{\sigma^2} | \nabla \sigma |^2
   = \io \frac1{\sigma^2} | \nabla \sigma |^2
   = \| \nabla \ln \sigma \|^2.
\end{equation}
On the other hand, 
\begin{equation}\label{co:25}
  \chi \io \frac{\alpha(\sigma)}{\sigma^2} \nabla\fhi \cdot \nabla \sigma
   = \chi \io \frac{1}{\sigma + \sigma^{p}} \nabla\fhi \cdot \nabla \sigma.
\end{equation}
Hence,
\begin{equation}\label{co:25b}
  \chi \io \frac{\alpha(\sigma)}{\sigma^2} \nabla\fhi \cdot \nabla \sigma
   \le \frac12 \io \frac1{\sigma^2} | \nabla \sigma |^2
    + c \| \nabla \fhi \|^2
    = \frac12 \| \nabla \ln\sigma \|^2
    + c \| \nabla \fhi \|^2.
\end{equation}
Replacing the above calculations into \eqref{co:24}, we deduce
\begin{equation}\label{co:24d}
  \ddt \io -\ln \sigma 
   + \frac12 \| \nabla \ln \sigma \|^2
  \le c \big( 1 + \| \nabla\fhi \|^2 \big).
\end{equation}
Summing \eqref{co:24d} to \eqref{en:14b} we then obtain
\begin{align}\no 
  & \ddt \Big(\frac12 \| \nabla \fhi \|^2 + \io F(\fhi) 
   + \io \big( \gammaciapo(\sigma) - \chi \sigma\fhi - \ln\sigma  \big) \Big)
   + \| \nabla\mu \|^2 
   + \io \alpha(\sigma) \big| \nabla ( \gamma(\sigma) - \chi \fhi ) \big|^2\\
 \label{en:14c}
  & \mbox{}~~~~~
  + \epsi \| D \vu \|^2 
   + \| \vu \|^2
   + \frac14 \| \nabla \ln \sigma \|^2
   \le c \big( 1 + \| \sigma \|_{L^p(\Omega)}^p 
     + \| \nabla\fhi \|^2 \big).
\end{align}
Then, recalling \eqref{gammaciapo} and \eqref{coer:energy}-\eqref{bdd:energy},
an application of Gr\"onwall's lemma 
permits us to deduce the following properties:
\begin{align}\label{st:11}
  & \| \fhi \|_{\LIV} \le c,\\
 \label{st:12}
  & \| F(\fhi) \|_{L^\infty(0,T;L^1(\Omega)} \le c,\\
 \label{st:13}
  & \| \nabla\mu \|_{\LDH} \le c,\\
 \label{st:14}
  & \| \vu \|_{\LDH} + \epsi^{1/2} \| D\vu \|_{\LDH} \le c,\\
 \label{st:15}
  & \| \sigma \|_{L^\infty(0,T;L^p(\Omega))} \le c,\\
 \label{st:cross}
  & \| \alpha^{1/2}(\sigma) \nabla (\gamma(\sigma) - \chi\fhi) \|_{\LDH} \le c,\\
 \label{st:26}
  & \| \ln \sigma \|_{L^\infty(0,T;L^1(\Omega))} + \| \ln \sigma \|_{\LDV} \le c.
\end{align}
Note that the above estimates, as well as the ones that follow, are to be intended
as a-priori bounds holding uniformly with respect to any hypothetical
regularization parameter. Moreover, as is implicitly written in \eqref{st:26}
we have in particular that $\sigma>0$ a.e.~in 
$Q$. Analogously, from \eqref{st:12} there {\sl formally}\/
follows that
\begin{equation}\label{st:a0}
   \| \fhi \|_{L^\infty(0,T;L^\infty(\Omega))} \le c.
\end{equation}
However, it is worth remarking that \eqref{st:a0} holds 
{\sl in the limit}, and with $c=1$, thanks to the specific
expression \eqref{Flog} of the logarithmic potential $F$,
which constrains $\fhi$ to take values into $[-1,1]$.
As $F$ is replaced by a regularizing family
$F_n$, \eqref{st:a0} in principle may be lost, or replaced by
a weaker information (cf.\ Remark~\ref{oss:limi1}). 
Nevertheless, in Section~\ref{sec:appro}, we shall prove that 
if $F_n$ are suitably designed, then the analogue of \eqref{st:a0}
keeps holding in the approximation, with the constant
$c$ (possibly larger than $1$), being independent of the
approximation parameter $n$. This fact, which has not been 
used in deducing the energy estimate (cf.~Remark~\ref{oss:limi1}),
will be crucial in the sequel.

That said, applying Korn's inequality, from \eqref{st:14} there follows that
\begin{equation}\label{st:a1}
   \epsi^{1/2} \| \vu \|_{L^2(0,T;V)} \le c.
\end{equation}
Next, we provide an estimate of the logarithmic term $f(\fhi)$.
To this aim, we test \eqref{CH2} by $\fhi-\fhi\OO$ so to obtain
\begin{equation}\label{co:11}
   \| \nabla \fhi \|^2
    + \io f(\fhi) (\fhi - \fhi\OO) 
    = \chi \io \sigma (\fhi - \fhi\OO) 
    + \io (\mu - \mu\OO) (\fhi - \fhi\OO).
\end{equation}
Then, accounting for \eqref{medie}, a well-known argument from \cite{MiZe} 
(which is robust with respect to regularizations of the potential)
permits us to see that 
\begin{equation}\label{co:11b}
  \io f(\fhi) (\fhi - \fhi\OO) 
   \ge \kappa \| f(\fhi) \|_{L^1(\Omega)} - c,
\end{equation}
where the constants $\kappa>0$ and $c\ge 0$ depend only on the expression 
of $f$ and on the parameter $\delta$ in \eqref{medie}. Then, using \eqref{co:11b}
and \eqref{st:a0} in \eqref{co:11}, and applying the 
Poincar\'e-Wirtinger inequality, it is not difficult to arrive at
\begin{equation}\label{co:12}
  \kappa \| f(\fhi) \|_{L^1(\Omega)}
    \le c \big( 1 + \| \sigma \|_{L^1(\Omega)}
    + \| \nabla \mu \| \big),
\end{equation}
whence, squaring, integrating in time, and recalling \eqref{st:13}
and \eqref{st:15}, we deduce
\begin{equation}\label{st:16}
  \| f(\fhi) \|_{L^2(0,T;L^1(\Omega))} \le c.
\end{equation}
Then, integrating \eqref{CH2} over $\Omega$ and 
using the above relation together with \eqref{st:15}, we obtain
\begin{equation}\label{st:16b}
  \| \mu\oo \|_{L^2(0,T)} \le c.
\end{equation}
Recalling \eqref{st:13}, we then have
\begin{equation}\label{st:17}
  \| \mu \|_{L^2(0,T;V)} \le c.
\end{equation}
We finally observe that, using \eqref{st:15} and \eqref{st:17}, we may 
reinterpret \eqref{CH2} as a time-dependent family of 
elliptic problems with maximal monotone nonlinearities. This fact
allows us the use of a nowadays standard 
regularity argument, which we just sketch for the reader's convenience.
Namely, one can test \eqref{CH2} by
$|\beta(\fhi)|^{p-1}\sign(\beta(\fhi))$, which is monotone with 
respect to $\fhi$ (recall that $\beta$ denotes the ``monotone
part'' of $f$). Then, using conditions \eqref{st:15} and
\eqref{st:17} (with Sobolev's embeddings), it is not difficult 
to arrive at
\begin{equation}\label{st:18}
  \| \fhi \|_{L^2(0,T;W^{2,p}(\Omega))} 
  + \| f(\fhi) \|_{L^2(0,T;L^p(\Omega))} 
   \le c.
\end{equation}
More precisely, the second bound is a direct consequence of the estimate,
while the first one follows by a further comparison of terms and 
applying $L^p$-elliptic regularity results of Agmon-Douglis-Nirenberg type. 
Note that this procedure will be repeated in 
the sequence with different (and better) exponents.
\beos 
One may wonder if, under stronger (compared to 
\eqref{hp:sigma0}) positivity conditions on the initial datum, 
relation \eqref{st:26} could be improved, getting for instance
a strictly positive uniform lower bound for $\sigma$. Actually,
a natural way to prove such a property would be to
test \eqref{nutr} by (a truncation of) 
$-\sigma^{-\eta}$ for larger values 
of the exponent $\eta$, and then possibly applying some 
iteration procedure. However, even though condition 
\eqref{hp:b} seems to allow such a procedure, the regularity 
\eqref{st:11} available at this level for $\nabla\fhi$ 
seems not sufficient in order to close the estimate.
\eddos


\section{Entropy-type estimates}
\label{sec:apriori2a}

The energy estimate derived in the previous section holds without any restriction
on the exponent $p$ and on the the dimension. However, if
$p$ is very close to $1$ (meaning that the cross diffusion term behaves 
at infinity similarly to the genuine Keller-Segel case), then,
at least for $d=3$, the outcome, in terms
of a-priori regularity, seems not sufficient in order to pass to the 
limit in the approximation. The main issue arises when
dealing with the product terms (more precisely, the transport term in \eqref{nutr} 
and the Korteweg term depending on $\sigma$ in \eqref{brink}), which may not
be identified under the sole regularity provided by 
the energy bound.

The following additional argument, based on the so-called ``entropy'' method 
for fourth order PDE's, is crucial in order to improve the available 
a-priori regularity and, in particular, to prove separate bounds for the 
summands in $\bH$ (cf.\ \eqref{rego:H}). We shall present two versions 
of the entropy bound,  corresponding to the regularity settings 
of Theorem~\ref{teo:main} and of Theorem~\ref{teo:main2}.
The first version, referring to the case $d=2$, is simpler (but likely
not optimal): indeed, in this situation for every $p>1$ we can ``directly''
obtain an estimate corresponding to the ``natural'' parabolic regularity
that would hold also in the ``linear'' case. We decided to start detailing
this argument because it gives a first idea of the method,
with reduced technicalities compared to the subsequent 3D proof.


\subsection{Entropy estimate: simple version for $d=2$}
\label{subsec:en1}

The basic strategy to obtain the entropy bound simply consists in testing
\eqref{CH2} by $-\Delta\fhi$. As done for the previous estimates,
we shall directly work on the original system \eqref{CH1}-\eqref{incompr}; 
in Section~\ref{sec:appro}, we will discuss how also this procedure might
be adapted so to fit the framework of a regularization scheme. 

That said, standard manipulations give
\begin{equation}\label{rig:11}
  \| \Delta\fhi \|^2 
   = - \io f'(\fhi) | \nabla\fhi |^2
   + \io \nabla\fhi \cdot \nabla \mu
    + \chi \io \nabla\fhi \cdot \nabla \sigma.
\end{equation}
Now, using \eqref{st:11} and recalling that $f'(\cdot) \ge -\lambda$,
we easily obtain
\begin{equation}\label{rig:12}
 - \io f'(\fhi) | \nabla\fhi |^2
   + \io \nabla\fhi \cdot \nabla \mu
    + \chi \io \nabla\fhi \cdot \nabla \sigma
    \le c \big( 1 + \| \nabla \mu \| + \| \nabla \sigma \| \big).
\end{equation}
Hence, replacing the above into \eqref{rig:11} and squaring, we deduce 
\begin{equation}\label{rig:13}
  \| \Delta\fhi \|^4 
   \le c_1 \big( 1 + \| \nabla \mu \|^2 + \| \nabla \sigma \|^2 \big),
\end{equation}
where the constant $c_1$ depends only on the outcome of the 
previous a-priori estimates (hence it is a computable quantity 
independent of the approximation parameters).

In order to control the last term on the \rhs, we test
\eqref{nutr} by $\sigma$. Then, assumption \eqref{hp:b} and simple
integrations by parts give
\begin{equation}\label{rig:21}
  \frac12 \ddt \| \sigma \|^2
   + \| \nabla\sigma \|^2
  \le c \big( 1 + \| \sigma \|^2 \big)
   + \chi \io \sigma^{2-p} |\nabla \sigma| | \nabla \fhi |.
\end{equation}
Now, the last term on the \rhs\ can be estimated as follows:
\begin{align}\no
  \chi \io \sigma^{2-p} | \nabla \sigma | | \nabla \fhi |
   & \le c \| \sigma^{1/2} \|_{L^{\frac{4(1+\epsilon)}{\epsilon}}(\Omega)}
    \| \sigma^{\frac{3-2p}2} \|_{L^{2(1+\epsilon)}(\Omega)}
    \| \nabla \sigma \|
    \| \nabla \fhi\|_{L^{\frac{4(1+\epsilon)}{\epsilon}}(\Omega)}\\
 \no 
   & \le c \| \sigma \|^{1/2}_{L^{\frac{2(1+\epsilon)}{\epsilon}}(\Omega)}
    \| \sigma^{\frac{3-2p}2} \|_{L^{2(1+\epsilon)}(\Omega)}
    \| \nabla \sigma \|
    \| \nabla \fhi\|^{\frac{\epsilon}{2(1+\epsilon)}}
    \| \fhi\|_{H^2(\Omega)}^{\frac{2+\epsilon}{2(1+\epsilon)}}\\
 \no 
   & \le c_\epsilon \| \sigma \|_V^{3/2}
    \| \sigma^{\frac{3-2p}2} \|_{L^{2(1+\epsilon)}(\Omega)}
    \Big( 1 + \| \Delta \fhi\|^{\frac{2+\epsilon}{2(1+\epsilon)}} \Big)\\
 \no
   & \le \frac12 \| \sigma \|_{V}^2
    + c_\epsilon \| \sigma^{\frac{3-2p}2} \|_{L^{2(1+\epsilon)}(\Omega)}^4
    \Big( 1 + \| \Delta\fhi\|^{\frac{2(2+\epsilon)}{1+\epsilon}} \Big)\\
 \label{rig:22-2d}
   & \le \frac12 \| \sigma \|_{V}^2
    + \delta \| \Delta \fhi\|^4
    + c 
    + c_{\epsilon,\delta} 
       \| \sigma^{\frac{3-2p}2} \|_{L^{2(1+\epsilon)}(\Omega)}^{\frac{8(1+\epsilon)}\epsilon}.
\end{align}
Here, we have used two dimensional embeddings with the Gagliardo-Nirenberg 
and Young inequalities. Moreover, $\epsilon>0$ and $\delta>0$ are ``small'' parameters
to be specified below, whereas, for instance, $c_\epsilon>0$ denotes a 
``large'' constant depending on the choice of $\epsilon$. To specify these
quantities, we first observe that, as $p>1$, we may use \eqref{st:15} to
control the last term on the \rhs, provided that $\epsilon>0$ is 
taken so small that 
\begin{align} \label{easy:epsi2}
  p \ge \frac{3(1+\epsilon)}{3+2\epsilon}, \quext{i.e.\ }\,
  (3-2p)(1+\epsilon)\le p.
\end{align}
Next, let us observe that, by elliptic regularity and 
\eqref{st:11}, 
\begin{equation}\label{rig:23a}
  \| \fhi \|_{H^2(\Omega)}^4 
   \le c_\Omega \big( \| \fhi \|^4 + \| \Delta\fhi \|^4 \big)
   \le c_\Omega \big( 1 + \| \Delta\fhi \|^4 \big).
\end{equation}
Then, replacing \eqref{rig:22-2d} into \eqref{rig:21}
and using \eqref{easy:epsi2}, we obtain
\begin{equation}\label{rig:23}
  \frac12 \ddt \| \sigma \|^2
   + \frac12 \| \nabla\sigma \|^2
  \le c \big( 1 + \| \sigma \|^2 \big)
   + c_\Omega \delta \| \Delta\fhi \|^4 
  + c_\delta \| \sigma \|_{L^{p}(\Omega)}^{{\boldsymbol P}},
\end{equation}
where ${\boldsymbol P}=4(1+\epsilon)(3-2p)/\epsilon$
is a ``large'', but otherwise fixed, exponent depending on $\epsilon$.
Next, recalling \eqref{rig:13}, we take $\delta>0$ 
so small that $2c_1 c_\Omega \delta \le 1/4$. Then, 
we multiply \eqref{rig:13} by $2c_\Omega\delta$
and subsequently sum the resulting relation to \eqref{rig:23} so
to obtain
\begin{equation}\label{rig:24}
  \frac12 \ddt \| \sigma \|^2
   + c_\Omega \delta \| \Delta\fhi \|^4 
   + \frac14\| \nabla\sigma \|^2
  \le c \big( 1 + \| \sigma \|^2 + \| \nabla \mu \|^2\big)
      + c_\delta \| \sigma \|_{L^{p}(\Omega)}^{{\boldsymbol P}}.
\end{equation}
Then, integrating in time and using Gr\"onwall's lemma with \eqref{st:17}
and \eqref{st:15}, we readily deduce the additional estimates
\begin{align}\label{st:ent11}
  & \| \fhi \|_{L^4(0,T;H^2(\Omega))} \le c,\\
 \label{st:ent12}
  & \| \sigma \|_{L^\infty(0,T;H)}
    + \| \sigma \|_{\LDV} \le c,
\end{align}
where we also used assumption \eqref{sigma0:en} on the initial datum. With the 
above regularity at hand, we can also improve the control of the logarithmic
term in \eqref{CH2}. In particular, using that, for $d=2$, $V\subset L^P(\Omega)$
for every $P\in[1,\infty)$, estimates \eqref{st:17} and \eqref{st:ent12}
and the same elliptic regularity argument used at the end of Section~\ref{sec:energy}
permit us to deduce 
\begin{equation}\label{st:ent2d}
  \| f(\fhi) \|_{L^2(0,T;L^P(\Omega))} +\| \fhi \|_{L^2(0,T;W^{2,P}(\Omega))}\le c_P%
    ~~\text{for all }\,P\in[1,\infty).
\end{equation}


\subsection{Entropy estimate: improved version for $d=3$}
\label{subsec:en2}

In this part we adapt the above argument so to deal with the assumptions
of Theorem~\ref{teo:main2} regarding the three-dimensional case.
In this setting, we will obtain a more accurate estimation of the cross-diffusion 
terms by valorizing the following chain of inequalities:
\begin{equation}\label{gani:ell}
  \| \nabla\fhi \|_{L^4(\Omega)} 
   \le c \| \fhi \|_{L^\infty(\Omega)}^{1/2} \| \fhi \|_{H^2(\Omega)}^{1/2}
   \le c \big(1 + \| \Delta\fhi \|^{1/2} \big),
\end{equation}
which holds in any space dimension and follows from
the Gagliardo-Nirenberg interpolation inequality \cite{Nir} 
and standard elliptic regularity results (recall that
$\fhi$ satisfies a no-flux boundary condition). Relation 
\eqref{gani:ell} uses in an essential way the uniform boundedness
of $\fhi$. This surely holds ``in the limit'' thanks to the structure
of the potential \eqref{Flog}, but it can be assumed to
be satisfied (uniformly) in the approximation as 
shown in Lemma~\ref{lem:linfty} below.

That said, we start, as before, by testing \eqref{CH2} by $-\Delta\fhi$
so to obtain the analogue of \eqref{rig:11}, namely
\begin{equation}\label{co:51}
   \| \Delta\fhi \|^2 
   \le c \big( 1 + \| \nabla\mu \| \big)
    + \chi \io \nabla \fhi \cdot \nabla \sigma.
\end{equation}
Moreover, we test \eqref{nutr2} by $\sigma^{p-1}$, with the purpose 
of obtaining a decoupled estimate for the cross-diffusion term.
We deduce
\begin{equation}\label{co:52}
  \frac1{p} \ddt \io \sigma^{p}
   + \frac{4(p-1)}{p^2} \io | \nabla \sigma^{\frac{p}2} |^2
   \le c + c \io \sigma^{p}
   + (p-1) \chi \io \frac{\sigma^{p-1}}{1+\sigma^{p-1}} \nabla\fhi\cdot\nabla \sigma.
\end{equation}
In order to control the terms on the \rhs\ of \eqref{co:51} and \eqref{co:52},
which are of the same type, we may proceed as follows:
\begin{align}\no
  & \io \nabla \fhi \cdot \nabla \sigma
    \le c \io | \nabla \fhi| | \nabla \sigma |
     \le c \io \sigma^{\frac{2-p}{2}} |\nabla \sigma^{p/2} | |\nabla\fhi|\\
  \no
  & \mbox{}~~~~~    
   \le c \| \sigma^{\frac{2-p}{2}} \|_{L^4(\Omega)} 
    \| \nabla \sigma^{p/2}  \| 
    \| \nabla \fhi \|_{L^4(\Omega)} \\
  \no
  & \mbox{}~~~~~    
   \le c \| \sigma^{\frac{2-p}{2}} \|_{L^{\frac{2p}{2-p}}(\Omega)}^{\eta}
     \| \sigma^{\frac{2-p}{2}} \|_{L^{\frac{6p}{2-p}}(\Omega)}^{1-\eta}
     \| \nabla \sigma^{p/2}  \| 
    \big( 1 + \| \Delta \fhi \|^{1/2} \big) \\
 \no
  & \mbox{}~~~~~    
   \le c \| \sigma \|_{L^p(\Omega)}^{\frac{\eta(2-p)}{2}}
     \| \sigma^{{\frac{p}2}{\frac{2-p}{p}}} \|_{L^{\frac{6p}{2-p}}(\Omega)}^{1-\eta}
     \| \nabla \sigma^{p/2}  \| 
    \big( 1 + \| \Delta \fhi \|^{1/2} \big) \\
 \label{co:54a}
  & \mbox{}~~~~~    
   \le c \| \sigma \|_{L^p(\Omega)}^{\frac{\eta(2-p)}{2}}
     \| \sigma^{{\frac{p}2}} \|_{L^6(\Omega)}^{\frac{(2-p)(1-\eta)}p}
     \| \nabla \sigma^{p/2}  \| 
    \big( 1 + \| \Delta \fhi \|^{1/2} \big),
\end{align}
where the interpolation exponent $\eta$ is
given by the relation
\begin{equation}\label{co:alpha}
  \eta\frac{2-p}{2p} + (1-\eta)\frac{2-p}{6p} = \frac14,
    \qquext{i.e.,\ }\,\eta=\frac{5p-4}{8-4p},
    \quad 1-\eta=\frac{12-9p}{8-4p}.
\end{equation}
Note that here we have implicitly assumed that $p\le 4/3$,
whence $2p/(2-p)\le 4$, which justifies the use of interpolation.
However, it is clear that, for $p>4/3$, the procedure can be modified 
as follows:
\begin{align}\no
  & \io \nabla \fhi \cdot \nabla \sigma
    \le c \| \sigma^{\frac{2-p}{2}} \|_{L^4(\Omega)} 
    \| \nabla \sigma^{p/2}  \| 
    \| \nabla \fhi \|_{L^4(\Omega)} \\
  \no
 & \mbox{}~~~~~    
  \le c \big( 1 + \| \sigma^{1/3} \|_{L^4(\Omega)} \big)
    \| \nabla \sigma^{p/2}  \| 
    \| \nabla \fhi \|_{L^4(\Omega)}\\
 \label{co:54b43}   
 & \mbox{}~~~~~  
  \le c \big( 1 + \| \sigma \|_{L^p(\Omega)}^{1/3} \big)
    \| \nabla \sigma^{p/2}  \| 
    \| \nabla \fhi \|_{L^4(\Omega)} 
  \le c \| \nabla \sigma^{p/2}  \| 
   \big( 1 + \| \Delta \fhi \|^{1/2} \big),
\end{align}  
the last inequality following from \eqref{st:15}. It is then not difficult
to adapt the remainder of the procedure so to cover also this (simpler)
situation.

Thus, going back to the case $p\le 4/3$, \eqref{co:54a} can be 
continuated as follows:
\begin{align}\no
  & \io \nabla \fhi \cdot \nabla \sigma
      \le c \| \sigma \|_{L^p(\Omega)}^{\frac{\eta(2-p)}{2}}
     \| \sigma^{{\frac{p}2}} \|_{L^6(\Omega)}^{\frac{(2-p)(1-\eta)}p}
     \| \nabla \sigma^{p/2}  \| 
    \big( 1 + \| \Delta \fhi \|^{1/2} \big)\\
 \label{co:54b}
  & \mbox{}~~~~~    
   \le c \| \sigma \|_{L^p(\Omega)}^{\frac{5p-4}{8}}
     \big( 1 + \| \nabla \sigma^{p/2} \|^{\frac{12-5p}{4p}} \big)
     \big( 1 + \| \Delta \fhi \|^{1/2} \big).
\end{align}
Replacing the above into \eqref{co:51}, we deduce 
\begin{equation}\label{co:51a} 
   \| \Delta\fhi \|^{3/2}
   \le c + c \| \nabla\mu \|^{3/4}
    + c \| \sigma \|_{L^p(\Omega)}^{\frac{5p-4}{8}} 
     \big( 1 + \| \nabla \sigma^{p/2} \|^{\frac{12-5p}{4p}} \big).
\end{equation}
Now, let us go back to \eqref{co:52}. Estimating the last term on 
the \rhs\ as done in \eqref{co:54b}, we have
\begin{equation}\label{co:52z}
  \frac1{p} \ddt \| \sigma \|_{L^p(\Omega)}^p
   + \frac{4(p-1)}{p^2} \| \nabla \sigma^{\frac{p}2} \|^2
   \le c + c \| \sigma \|_{L^p(\Omega)}^p
     + c \| \sigma \|_{L^p(\Omega)}^{\frac{5p-4}{8}}
     \big( 1 + \| \nabla \sigma^{p/2} \|^{\frac{12-5p}{4p}} \big)
     \big( 1 + \| \Delta \fhi \|^{1/2} \big).
\end{equation}
Using \eqref{co:51a} to control the term with the Laplacian,
we then infer 
\begin{align}\no
  & \frac1{p} \ddt \| \sigma \|_{L^p(\Omega)}^p
   + \frac{4(p-1)}{p^2} \| \nabla \sigma^{\frac{p}2} \|^2
   \le c + \bC\\
 \label{co:x1}
  & \mbox{}~~~~~~~~~~
     + \bC \big( 1 + \| \nabla \sigma^{p/2} \|^{\frac{12-5p}{4p}} \big)
     \Big( 1 + \| \nabla\mu \|^{1/4}
    + \bC \big( 1 + \| \nabla \sigma^{p/2} \|^{\frac{12-5p}{12p}} \big) \Big).
\end{align}
Here and below, $\bC$ represents a ``large'', but otherwise computable,
constant depending on (some power of) the $L^p(\Omega)$-norm of $\sigma$,
which has already been estimated uniformly in time thanks to \eqref{st:15}.

Now, we first observe that, as $p>12/11$, it follows 
\begin{equation}\label{co:x2}
  \frac{12-5p}{4p} < \frac74, \quext{whence }~
  \bC \| \nabla \sigma^{p/2} \|^{\frac{12-5p}{4p}} \| \nabla\mu \|^{1/4}
   \le \frac{p-1}{p^2}\| \nabla \sigma^{p/2} \|^2 + \bC \| \nabla\mu \|^2
   + c.
\end{equation}
Moreover, the condition $p>12/11$ also guarantees
\begin{equation}\label{co:x3x}
  \frac{12-5p}{4p} + \frac{12-5p}{12p} = \frac{12-5p}{3p}<2, 
  \quext{whence }~
  \bC \| \nabla \sigma^{p/2} \|^{\frac{12-5p}{3p}} 
   \le \frac{p-1}{p^2}\| \nabla \sigma^{p/2} \|^2 + \bC.
\end{equation}
On account of these considerations, \eqref{co:x1} gives
\begin{equation}\label{co:x4}
 \frac1{p} \ddt \| \sigma \|_{L^p(\Omega)}^p
   + \frac{2(p-1)}{p^2} \| \nabla \sigma^{\frac{p}2} \|^2
   \le \bC \big( 1 + \| \nabla\mu \|^2 \big).
\end{equation}
Hence, applying Gr\"onwall's lemma, we obtain
\begin{equation}\label{st:ent22}
  \| \sigma^{p/2} \|_{\LDV} \le c,
\end{equation}
which implies the desired decoupled regularity bound for the cross-diffusion
term in \eqref{nutr}. Finally, going back to \eqref{co:51a}, we prove an estimate 
of $\Delta\fhi$. To this aim, we need to distinguish
between two case: if $p\le 3/2$, then taking a suitable power of that relation
we deduce
\begin{equation}\label{co:51b} 
   \| \Delta\fhi \|^{\frac{12p}{12-5p}}
   \le c + c \| \nabla\mu \|^{\frac{6p}{12-5p}}
    + c \| \sigma \|_{L^p(\Omega)}^{\frac{p(5p-4)}{12-5p}}
     \big( 1 + \| \nabla \sigma^{p/2} \|^2 \big),
\end{equation}
and it turns out that $\frac{6p}{12-5p}\le 2$. Hence, integrating in
time and using \eqref{st:ent22} we infer
\begin{equation}\label{st:ent21}
  \| \fhi \|_{L^{\frac{12p}{12-5p}}(0,T;H^2(\Omega))} \le c.
\end{equation}
On the other hand, if $p> 3/2$, we take the $(8/3)$-power of \eqref{co:51a},
which gives 
\begin{equation}\label{co:51bb} 
   \| \Delta\fhi \|^{4}
   \le c + c \| \nabla\mu \|^2
    + c \| \sigma \|_{L^p(\Omega)}^{\frac{5p-4}{3}}
     \big( 1 + \| \nabla \sigma^{p/2} \|^{\frac{24-10p}{3p}} \big)
\end{equation}
Then, $p>3/2$ implies $\frac{24-10p}{3p} <2$. Hence, integrating once more
in time and using \eqref{st:ent22}, in place of \eqref{st:ent21} we get 
\begin{equation}\label{st:ent21b}
  \| \fhi \|_{L^4(0,T;H^2(\Omega))} \le c.
\end{equation}


\subsection{Additional regularity for $d=3$ by iteration of the entropy estimate}
\label{subsec:boot}

We now see how the previous procedure can be adapted in order to achieve
additional summability of $\sigma$ in the case when 
$\sigma_0\in L^{q}(\Omega)$. For technical reasons, we will start
assuming $p<q\le 20/11$. For larger $q$ the procedure will be 
carried out by means of a two-step iteration.
That said, we start with testing \eqref{nutr2} 
by $\sigma^{q-1}$. Then, using once more \eqref{hp:b}, we deduce
\begin{align}\no
  \frac1{q} \ddt \io \sigma^{q}
   + \frac{4(q-1)}{q^2} \io | \nabla \sigma^{\frac{q}2} |^2
  & \le c + c \io \sigma^{q}
   + (q-1) \chi \io \frac{\sigma^{q-1}}{1+\sigma^{p-1}} \nabla\fhi\cdot\nabla \sigma\\
 \label{co:52q2}  
  & \le c + c \io \sigma^{q}
   + c \io \sigma^{q-p} |\nabla\fhi| |\nabla \sigma|.
\end{align}
In order to estimate the last term on the \rhs\ in an optimal way,
we go back to relation \eqref{co:51} and modify \eqref{co:54a}
as follows:
\begin{align}\no
  & \chi \io \nabla \fhi \cdot \nabla \sigma
     \le c \io \sigma^{\frac{2-q}{2}} |\nabla \sigma^{q/2} | |\nabla\fhi|\\
  \no
  & \mbox{}~~~~~    
   \le c \| \sigma^{\frac{2-q}{2}} \|_{L^4(\Omega)} 
    \| \nabla \sigma^{q/2}  \| 
    \| \nabla \fhi \|_{L^4(\Omega)} \\
  \no
  & \mbox{}~~~~~    
   \le c \| \sigma^{\frac{2-q}{2}} \|_{L^{\frac{2p}{2-q}}(\Omega)}^{\eta}
     \| \sigma^{\frac{2-q}{2}} \|_{L^{\frac{6q}{2-q}}(\Omega)}^{1-\eta}
     \| \nabla \sigma^{q/2}  \| 
    \big( 1 + \| \Delta \fhi \|^{1/2} \big) \\
 \no
  & \mbox{}~~~~~    
   \le c \| \sigma \|_{L^p(\Omega)}^{\frac{\eta(2-q)}{2}}
     \| \sigma^{{\frac{q}2}{\frac{2-q}{q}}} \|_{L^{\frac{6q}{2-q}}(\Omega)}^{1-\eta}
     \| \nabla \sigma^{q/2}  \| 
    \big( 1 + \| \Delta \fhi \|^{1/2} \big) \\
 \label{co:54arb}
  & \mbox{}~~~~~    
   \le c \| \sigma \|_{L^p(\Omega)}^{\frac{\eta(2-q)}{2}}
     \| \sigma^{{\frac{q}2}} \|_{L^6(\Omega)}^{\frac{(2-q)(1-\eta)}q}
     \| \nabla \sigma^{q/2}  \| 
    \big( 1 + \| \Delta \fhi \|^{1/2} \big).
\end{align}
As before, the previous computation is rigorous provided that
\begin{equation}\label{cond:12}
  \frac{2p}{2-q}<4<\frac{6q}{2-q}.
\end{equation}
Clearly the right inequality is always satisfied as far as $1<q<2$. Concerning
the left one, it corresponds to asking $2-q>p/2$. On the other hand, 
if $2-q\le p/2$, as before we may notice that 
\begin{align}\no
  & \chi \io \nabla \fhi \cdot \nabla \sigma
     \le c \io \sigma^{\frac{2-q}{2}} |\nabla \sigma^{q/2} | |\nabla\fhi|\\
  \no
  & \mbox{}~~~~~    
   \le c \| \sigma^{\frac{2-q}{2}} \|_{L^4(\Omega)} 
    \| \nabla \sigma^{q/2}  \| 
    \| \nabla \fhi \|_{L^4(\Omega)} 
  \le c \big( 1 + \| \sigma^{\frac{p}4} \|_{L^4(\Omega)} \big)
    \| \nabla \sigma^{q/2}  \| 
    \| \nabla \fhi \|_{L^4(\Omega)} \\
 \label{co:54arb2}
  & \mbox{}~~~~~    
   \le \bC
    \| \nabla \sigma^{q/2}  \| 
   \big( 1 + \| \Delta \fhi \|^{1/2} \big),
\end{align}
and the argument proceeds with straighforward adaptations.

That said, going back to the case $\frac{2p}{2-q}<4$, 
the interpolation exponent $\eta$ is then provided by the relation
\begin{equation}\label{co:alpharb}
  \eta\frac{2-q}{2p} + (1-\eta)\frac{2-q}{6q} = \frac14,
    \qquext{i.e.,\ }\,\eta=\frac{(5q-4)p}{2(2-q)(3q-p)},
    \quad 1-\eta=\frac{q(12-6q-3p)}{2(2-q)(3q-p)}.
\end{equation}
Hence, \eqref{co:54arb} can be continuated as follows:
\begin{align}\no
  & \chi \io \nabla \fhi \cdot \nabla \sigma
      \le c \| \sigma \|_{L^p(\Omega)}^{\frac{\eta(2-q)}{2}}
     \| \sigma^{{\frac{q}2}} \|_{L^6(\Omega)}^{\frac{(2-q)(1-\eta)}q}
     \| \nabla \sigma^{q/2}  \| 
    \big( 1 + \| \Delta \fhi \|^{1/2} \big)\\
 \no
  & \mbox{}~~~~~    
   \le c \| \sigma \|_{L^p(\Omega)}^{\frac{(5q-4)p}{4(3q-p)}}
     \big( 1 + \| \nabla \sigma^{q/2} \|^{\frac{12-6q-3p}{6q-2p}} \big)
     \| \nabla \sigma^{q/2}  \| 
     \big( 1 + \| \Delta \fhi \|^{1/2} \big)\\
 \label{co:54brb}
  & \mbox{}~~~~~    
   \le \bC
     \big( 1 + \| \nabla \sigma^{q/2} \|^{\frac{12-5p}{6q-2p}} \big)
     \big( 1 + \| \Delta \fhi \|^{1/2} \big),
\end{align}
where, as before, $\bC>0$ is a (possibly large) constant 
depending on the outcome of the previous estimates (specifically, 
on \eqref{st:15}).

Hence, in view of the above procedure, \eqref{co:51a} is replaced by 
\begin{equation}\label{co:51abb}
   \| \Delta\fhi \|^{3/2}
   \le c + c \| \nabla\mu \|^{3/4}
    + \bC \big( 1 + \| \nabla \sigma^{q/2} \|^{\frac{12-5p}{6q-2p}} \big),
\end{equation}
or, in other words,
\begin{equation}\label{co:51abc}
   \| \Delta\fhi \|^{\frac{18q-6p}{12-5p}}
   \le c + c \| \nabla\mu \|^{\frac{9q-3p}{12-5p}}
    + \bC \big( 1 + \| \nabla \sigma^{q/2} \|^{2} \big).
\end{equation}
Now, let us go back to \eqref{co:52q2}. We propose the following estimation
strategy:
\begin{align}\no
  & c \io \sigma^{q-p} | \nabla \sigma | | \nabla \fhi |
   \le c \io \sigma^{\frac{q-2p+2}{2}} | \nabla \sigma^{q/2} | | \nabla \fhi |\\
 \no
  & \mbox{}~~~~~
   = c \io \sigma^{\frac{q}2\frac{q-2p+2}{q}} | \nabla \sigma^{q/2} | | \nabla \fhi |
   \le c \| \sigma^{\frac{q}2\frac{q-2p+2}{q}} \|_{L^4(\Omega)}
   \| \nabla \sigma^{q/2} \| \| \nabla \fhi \|_{L^4(\Omega)}\\
 \label{co:y11}
  & \mbox{}~~~~~
    \le c \| \sigma^{q/2} \|_{L^{\frac{4(q-2p+2)}{q}}(\Omega)}^{\frac{q-2p+2}q}
     \| \nabla \sigma^{q/2} \| \big( 1 +  \| \Delta \fhi \|^{1/2} \big).
\end{align}
The first factor on the \rhs\ is controlled by interpolation with the following 
choice:
\begin{equation}\label{co:y12}
   \frac{q}{4(q-2p+2)} = \eta \frac{q}{2p} + (1-\eta) \frac16,
\end{equation}
whence 
\begin{align}\label{co:y13}
  & \eta = \frac{pq+4p^2-4p}{(q-2p+2)(6q-2p)} = \frac{pq+4p^2-4p}{6q^2+4p^2-14pq-4p+12q}, \\
 \label{co:y14}
  & 1 - \eta = \frac{6q^2-15pq+12q}{(q-2p+2)(6q-2p)} = \frac{6q^2-15pq+12q}{6q^2+4p^2-14pq-4p+12q}.
\end{align}
To be precise, in deducing the last inequality in \eqref{co:y11} we have implicitly
used that $\frac{4(q-2p+2)}{q}\ge 1$ (which is always true), while in \eqref{co:y12} 
we have used the more restrictive condition
\begin{equation}\label{co:y12b}
   \frac{2p}{q}\le \frac{4(q-2p+2)}q \le 6.
\end{equation}
Actually, it is easy to check that the right inequality in \eqref{co:y12b} is always satisfied
in our exponents range. This is not always true for the left inequality, instead; nevertheless,
if the opposite holds, then it follows $2(q-2p+2)<p$. Hence, going back to \eqref{co:y11}, 
one has
\begin{align}\no
  & c \io \sigma^{q-p} | \nabla \sigma | | \nabla \fhi |
   \le c \| \sigma^{\frac{q-2p+2}{2}} \|_{L^4(\Omega)}
   \| \nabla \sigma^{q/2} \| \| \nabla \fhi \|_{L^4(\Omega)}\\
 \no
  & \mbox{}~~~~~
   =  c \bigg( \io \sigma^{2(q-2p+2)} \bigg)^{1/4}
   \| \nabla \sigma^{q/2} \| \| \nabla \fhi \|_{L^4(\Omega)}
  \le  c \big( 1 + \| \sigma \|_{L^p(\Omega)}^{p/4} \big)
   \| \nabla \sigma^{q/2} \| \| \nabla \fhi \|_{L^4(\Omega)}\\
 \label{co:y11x}
 & \mbox{}~~~~~
  \le  \bC \| \nabla \sigma^{q/2} \| \big( 1 +  \| \Delta \fhi \|^{1/2} \big),
\end{align}
and the remainder of the procedure can be adapted up to small variations.

Hence, going back to the (more delicate) case when \eqref{co:y12b} holds, 
\eqref{co:y11} can be continuated as follows:
\begin{align}\no
  & c \io \sigma^{q-p} | \nabla \sigma | | \nabla \fhi |\\
 \no
  & \mbox{}~~~~~
   \le c \| \sigma^{q/2} \|_{L^{\frac{2p}q}(\Omega)}^{\frac{\eta(q-2p+2)}q}
    \| \sigma^{q/2} \|_{L^6(\Omega)}^{(1-\eta)\frac{q-2p+2}q}
     \| \nabla \sigma^{q/2} \| \big( 1 +  \| \Delta \fhi \|^{1/2} \big)\\
 \no
  & \mbox{}~~~~~
   \le c \| \sigma^{q/2} \|_{L^{\frac{2p}q}(\Omega)}^{\frac{pq+4p^2-4p}{q(6q-2p)}}
    \| \sigma^{q/2} \|_V^{\frac{6q-15p+12}{6q-2p}}
     \| \nabla \sigma^{q/2} \| \big( 1 +  \| \Delta \fhi \|^{1/2} \big)\\
 \label{co:y16}
  & \mbox{}~~~~~
   \le c \| \sigma \|_{L^p(\Omega)}^{\frac{pq+4p^2-4p}{2(6q-2p)}}
    \big ( 1 + \| \nabla \sigma^{q/2} \|^{\frac{12q-17p+12}{6q-2p}} \big)
        \big( 1 +  \| \Delta \fhi \|^{1/2} \big).
\end{align}
At this point, we replace \eqref{co:y16} into \eqref{co:52q2}
and subsequently estimate the term with the Laplacian by means of 
\eqref{co:51abb}. We obtain
\begin{align}\no
  & \frac1{q} \ddt \| \sigma \|_{L^q(\Omega)}^q
   + \frac{4(q-1)}{q^2} \| \nabla \sigma^{\frac{q}2} \|^2
   \le c + c \| \sigma \|_{L^q(\Omega)}^q
   + \bC \big ( 1 + \| \nabla \sigma^{q/2} \|^{\frac{12q-17p+12}{6q-2p}} \big)
        \big( 1 +  \| \Delta \fhi \|^{1/2} \big)\\
 \label{co:y17}
  & \mbox{}~~~~~
   \le c + c \| \sigma \|_{L^q(\Omega)}^q
   + \bC \big ( 1 + \| \nabla \sigma^{q/2} \|^{\frac{12q-17p+12}{6q-2p}} \big)
        \big( 1 + \| \nabla\mu \|^{1/4}
        + \| \nabla \sigma^{q/2} \|^{\frac{12-5p}{18q-6p}} \big).
\end{align}
Now, in order to control the \rhs\ via Gr\"onwall's lemma we need two conditions. 
The first one, in analogy with \eqref{co:x2}, corresponds to 
\begin{equation}\label{co:x2b}
  \frac{12q-17p+12}{6q-2p} \le \frac74, \quext{or, equivalently, }\,
    q \le 9p - 8.
\end{equation}
It is then easy to check that the above inequality holds 
under the technical condition $p<q\le 20/11$ (actually, if $p\ge 10/9$ 
then it holds for every $q\in(p,2]$). This implies 
\begin{equation}\label{co:x2c}
  \bC \| \nabla \sigma^{q/2} \|^{\frac{12q-17p+12}{6q-2p}} 
      \| \nabla\mu \|^{1/4}
   \le \frac{q-1}{q^2}\| \nabla \sigma^{q/2} \|^2 + \bC \| \nabla\mu \|^2
    + c.
\end{equation}
Moreover, let us compute
\begin{equation}\label{co:x3}
  \frac{12q-17p+12}{6q-2p} + \frac{12-5p}{18q-6p}
   = \frac{18q-28p+24}{9q-3p}.
\end{equation}
Then, we notice that 
\begin{equation}\label{co:x3b}
  \frac{18q-28p+24}{9q-3p}<2 \quext{corresponds exactly to }\,
   22p > 24,~~\text{i.e.~}\,p>12/11.
\end{equation}
Hence, in our admissible range of exponents, independently of the 
value assumed by $q$, we have
\begin{equation}\label{co:y18}
  \bC \big ( 1 + \| \nabla \sigma^{q/2} \|^{\frac{12q-17p+12}{6q-2p}} \big)
        \big( 1 + \| \nabla \sigma^{q/2} \|^{\frac{12-5p}{18q-6p}} \big)
   \le \bC + \frac{q-1}{q^2}\| \nabla \sigma^{q/2} \|^2.
\end{equation}
Replacing \eqref{co:x2c} and \eqref{co:y18} into \eqref{co:y17},
we then arrive at 
\begin{equation}\label{co:y19}
  \frac1{q} \ddt \| \sigma \|_{L^q(\Omega)}^q
   + \frac{2(q-1)}{q^2} \| \nabla \sigma^{\frac{q}2} \|^2
   \le c + c \| \sigma \|_{L^q(\Omega)}^q
   + \bC \big( 1 + \| \nabla\mu \|^2 \big).
\end{equation}
Applying Gr\"onwall's lemma and recalling \eqref{sigma0:en2}, we obtain 
\begin{equation}\label{st:ent22q}
  \| \sigma^{q/2} \|_{\LIH} + \| \sigma^{q/2} \|_{\LDV} \le c.
\end{equation}
To conclude, we deduce some additional regularity for $\Delta\fhi$. To this aim,
we go back to \eqref{co:51abc} and notice that, if $\frac{18q-6p}{12-5p}< 4$
then the exponent of the norm of $\nabla\mu$ on the \rhs\ is less than
$2$. Integrating in time, we then deduce 
\begin{equation}\label{st:ent21q}
  \| \fhi \|_{L^{\frac{18q-6p}{12-5p}}(0,T;H^2(\Omega))} \le c.
\end{equation}
On the other hand, if $\frac{18q-6p}{12-5p} \ge 4$, as before we need to rescale the 
inequality before integrating in time. In that case, we obtain 
\begin{equation}\label{st:ent21q4}
  \| \fhi \|_{L^{4}(0,T;H^2(\Omega))} \le c.
\end{equation}
Notice that the above conditions correspond exactly to the choice of the 
exponent $P_0$ in the statement of Theorem~\ref{teo:main2} 
(cf.\ in particular the first \eqref{P0}).

\smallskip

Let us now remove the technical condition $p<q\le 20/11$ and, for simplicity,
let us just consider the worst case scenario, i.e., when $p\sim 12/11$. 
In that case, let us assume, as stated in \eqref{sigma0:en2}, that 
$\sigma_0\in L^q(\Omega)$ for $q>20/11$. In this case, we proceed by
bootstrap: namely, we first perform the above argument for 
$q_0=20/11$ (we note this ``temporary'' exponent as $q_0$ so to
distinguish it from $q$ which now corresponds to the summability
of the initial datum). Then, we achieve the regularity properties
\eqref{st:ent22q} (where $q_0=20/11$) and \eqref{st:ent21q4}
(note that, as $q_0=20/11$, then $\frac{18q_0-6p}{12-5p} \ge 4$ for any 
$p\in (12/11,q_0)$). Then, to bootstrap regularity, assuming for simplicity
$q = 2$ (the ``intermediate'' case when $\sigma_0$ lies in $L^q(\Omega)$ 
for $q\in(20/11,2)$ can be treated with small modifications) one 
may achieve the ``linear'' parabolic regularity simply by 
testing \eqref{nutr} by $\sigma$. This gives
\begin{equation}\label{co:y20}
  \frac12 \ddt \| \sigma \|^2
   + \| \nabla \sigma \|^2
  \le c\big( 1 + \| \sigma \|^2 \big)
   + c \io \sigma^{2-p} \nabla \sigma\cdot \nabla \fhi.
\end{equation}
The last term, being $p>12/11$, can be controlled this way:
\begin{align}\no
  & \io \sigma^{2-p} \nabla \sigma\cdot \nabla \fhi
   \le c \big( 1 + \| \sigma^{10/11} \|_{L^4(\Omega)} \big) \| \nabla \sigma \| 
    \| \nabla \fhi \|_{L^4(\Omega)}\\
 \no
  & \mbox{}~~~~~
   \le c \big( 1 + \| \sigma^{10/11} \|_{L^6(\Omega)}^{3/4} \| \sigma^{10/11} \|^{1/4} \big)
    \| \nabla \sigma \| 
    \big( 1 + \| \Delta \fhi \|^{1/2} \big)\\
 \label{co:y32}
  & \mbox{}~~~~~
   \le c \big( 1 + \| \sigma \|_{L^6(\Omega)}^{3/4} \big)
    \| \nabla \sigma \| \big( 1 + \| \Delta \fhi \|^{1/2} \big)
   \le c + c \| \Delta \fhi \|^4 + \frac12 \| \nabla \sigma \|^2 .
\end{align}
Hence, integrating in time and using \eqref{sigma0:en2}, we obtain once 
more \eqref{st:ent22q} (for $q=2$). With simple adaptations,
one can get the analogue also in the case $q\in(20/11,2)$. 
By interpolation, \eqref{st:ent22q} implies
\begin{equation}\label{st:q1}
  \| \sigma \|_{L^{2}(0,T;L^{\frac{3q}{3-q}}(\Omega))}
   \le c.   
\end{equation} 
Hence, interpreting once more \eqref{CH2} as a time-dependent family of 
elliptic problems and using \eqref{st:q1} together with \eqref{st:17},
the usual regularity argument gives
\begin{equation}\label{st:y32}
  \| \fhi \|_{L^2(0,T;W^{2,\frac{3q}{3-q}}(\Omega))} 
   + \| f(\fhi) \|_{L^2(0,T;L^{\frac{3q}{3-q}}(\Omega))} \le c.
\end{equation}
%


\section{Passage to the limit}
\label{sec:lim}

In this section we complete the proof of our results. We will 
proceed by showing the so-called ``weak sequential
stability'' of solution families. Namely, we let $(\fhi_n,\mu_n,\sigma_n,\vu_n)$,
$n\in\NN$, be a sequence of solutions complying with the a-priori estimates
obtained before uniformly with respect to $n$. Then, proving weak sequential 
stability means that at least a (non-relabelled) subsequence has to
converge to a quadruple $(\fhi,\mu,\sigma,\vu)$ solving \eqref{CH1}-\eqref{incompr}
in the sense specified in the statements of the existence theorems.

Of course, in order for the argument to be fully rigorous, we should
rather assume that the sequence $(\fhi_n,\mu_n,\sigma_n,\vu_n)$ 
solves some approximation of the system (for instance the one 
sketched in Section~\ref{sec:appro} below), with the approximating
parameters varying with $n$ in a suitable way. However, for the sake of 
simplicity, we prefer to assume that $(\fhi_n,\mu_n,\sigma_n,\vu_n)$ complies 
directly with the ``limit'', or ``original'', system, possibly written 
in the ``strong'' form \eqref{CH1}-\eqref{incompr}. We actually believe
that such a simplified method permits us to focus on the 
substantial points of the compactness argument. In the next 
section, we will explain how the procedure could be adapted in order to
fit the proposed approximation scheme.

We will mostly focus on the (more difficult) case of Theorem~\ref{teo:main2},
referring to the situation where $d=3$ and $p>12/11$ (at the end we
will sketch the differences occurring for $d=2$ in the setting
of Theorem~\ref{teo:main}). Moreover, 
in order to achieve the limit $n\nearrow\infty$, we will use only the 
information resulting from the mass conservation, the minimum 
principle, the energy estimate, and the most general version of 
the entropy bound, namely (cf.\ \eqref{st:ent21} and \eqref{st:ent22})
\begin{equation}\label{st:n1}
   \| \sigma_n^{p/2} \|_{\LDV} 
    + \| \fhi_n \|_{L^{\frac{12p}{12-5p}}(0,T;H^2(\Omega))} \le c,
\end{equation}
where the exponent $\frac{12p}{12-5p}$ is replaced by $4$ for $p\ge 3/2$,
cf.~\eqref{st:ent21b}). Indeed, the additional information resulting 
from more refined versions of the entropy estimate, and in particular 
the argument performed in Subsec.~\ref{subsec:boot}, 
will produce additional regularity of weak solution, but is not 
essential for taking the limit. Note also that, from now on, we 
will stress the dependence on $n$ in the estimates. Moreover,
in order to prepare the asymptotic limit $\epsi\searrow 0$, 
we will note by $c$ (respectively $c\ee$) the
positive constants occurring in the estimates which are
independent of $\epsi$ (respectively,
depending on $\epsi$). 

That said, we start testing \eqref{CH1} by $\xi\in V$ and integrating by
parts so to obtain the weak formulation \eqref{CH1w}. Next,
by \eqref{st:14} and the uniform boundedness of $\fhi_n$ there follows
\begin{equation}\label{st:a2}
   \| \fhi_n \vu_n \|_{L^2(0,T;H)} \le c.
\end{equation}
Hence, comparing terms in \eqref{CH1}, 
and using \eqref{st:17}, \eqref{hp:h} and \eqref{st:a2}, 
we obtain 
\begin{equation}\label{st:a3}
   \| \fhi_{n,t} \|_{L^2(0,T;V')} \le c.
\end{equation}
Let us now test \eqref{nutr} by $\eta\in W^{1,\infty}(\Omega)$ 
(actually, in the procedure it will be clear that the conditions on 
$\eta$ can be weakened, as specified in the statements)
and integrate by parts so to deduce the weak formulation \eqref{nutrw}.
Recall in particular that we have set
\begin{equation}\label{defi:H}
  \bH(\sigma_n,\fhi_n) := \alpha^{1/2}(\sigma_n) \nabla( \gamma(\sigma_n) - \chi \fhi_n ),
\end{equation}
and, with this notation, \eqref{st:cross} corresponds to 
\begin{equation}\label{st:a4}
   \| \bH(\sigma_n,\fhi_n) \|_{L^2(0,T;H)} \le c.
\end{equation}
Then, observing that (we assume with no loss of generality $p\in(1,2)$,
the analogue being trivial in the case $p=2$)
\begin{equation}\label{st:a5}
  \| \alpha^{1/2}(\sigma_n) \|_{L^\infty(0,T;L^{\frac{2p}{2-p}}(\Omega))} 
   \le \| \sigma_n^{\frac{2-p}2} \|_{L^\infty(0,T;L^{\frac{2p}{2-p}}(\Omega))}    
   \le \| \sigma_n \|_{L^\infty(0,T;L^p(\Omega))}^{\frac{2-p}2}    
   \le c,  
\end{equation}
combining \eqref{st:a4} and \eqref{st:a5}, it is easy to deduce 
\begin{equation}\label{st:a6}
   \| \alpha^{1/2}(\sigma_n) \bH(\sigma_n,\fhi_n) \|_{L^2(0,T;L^p(\Omega))} \le c.
\end{equation}
Note that the above only relies on the ``energy'' estimate: using the ``entropy''
regularity, we could actually get some better information.

Next, we observe that \eqref{st:15} and the first \eqref{st:n1},
thanks to 
%
%
interpolation and Sobolev's embeddings, imply
\begin{equation}\label{st:n3}
  \| \sigma_n \|_{L^{\frac{5p}{3}}(Q)} 
  + \| \sigma_n \|_{L^{\frac{4p}{6-3p}}(0,T;H)}
  + \| \sigma_n \|_{L^{2}(0,T;L^{\frac{3p}{3-p}}(\Omega))}
   \le c.   
\end{equation} 
In particular, we may notice that the last two exponents depending 
on $p$ are strictly larger than $2$ if and only if $p$ is strictly 
larger than $6/5$.

Next, we observe that, from \eqref{st:11}, the second \eqref{st:n1} 
(note that $\frac{12p}{12-5p}>2$ for $p>12/11$)
and interpolation, there follows at least 
\begin{equation}\label{st:fhi:3d}
  \| \nabla \fhi_n \|_{L^{10/3}(Q)} 
   \le c \big( \| \fhi_n \|_{\LIV} + \| \fhi_n \|_{\LDHD} \big)
   \le c.   
\end{equation}
Moreover, using again \eqref{st:n1} with the first \eqref{st:a5}, we get 
\begin{equation}\label{st:fhi:3d2}
  \| \alpha(\sigma_n) \nabla \fhi_n \|_{L^2(0,T;L^{\frac{6p}{12-5p}}(\Omega))} 
    \le c \| \alpha(\sigma_n) \|_{L^\infty(0,T;L^{\frac{p}{2-p}}(\Omega))}
     \| \nabla\fhi_n \|_{L^2(0,T;L^6(\Omega))} 
   \le c.   
\end{equation}
Since $\frac{6p}{12-5p}>1$ for $p>12/11$, the above, combined with 
\eqref{st:a6}, implies
\begin{align}\no
  \| \nabla \sigma_n \|_{L^2(0,T;L^{S}(\Omega))} 
    & \le c \big( \| \alpha^{1/2}(\sigma_n) \bH(\sigma_n,\fhi_n) \|_{L^2(0,T;L^p(\Omega))}
    + \| \alpha(\sigma_n) \nabla \fhi_n \|_{L^2(0,T;L^{\frac{6p}{12-5p}}(\Omega))} \big)\\
  \label{st:fhi:3d3}
    & \le c, \qquext{where }\ S = \min\Big\{ \frac{6p}{12-5p}, p \Big\}>1.
\end{align}
We now consider the convection term and note that its control is the only step 
where the condition $p>12/11$ is not sufficient to obtain an estimate which
also allows for the Brinkman-to-Darcy limit $\epsi\searrow0$. Actually,
as $\epsi>0$ is fixed, we may use \eqref{st:a1},
which, combined with the {\sl third}\/ \eqref{st:n3}, thanks to Sobolev's embeddings,
implies
\begin{equation}\label{st:d12}
  \| \sigma_n\vu_n \|_{L^{Z}(0,T;L^{\frac43}(\Omega))} 
   \le c_\epsi,
\end{equation}
where $Z$ is some exponent {\sl strictly}\/ larger than $1$. To deduce the
above we used that, for $p=12/11$, one has 
$\frac{3p}{3-p}=12/7$, so that $1/6+7/12=3/4$. Since in fact $p>12/11$,
changing a bit the interpolation exponents in
\eqref{st:n3}, we see that $Z$ can be taken strictly larger than $1$
in \eqref{st:d12}. 

On the other hand, as we deal with the Darcy limit, in order to use the {\sl first}\/
\eqref{st:14}, we need to know that 
\begin{equation}\label{st:d13}
  \| \sigma_n \|_{L^{Z_1}(Q)} \le c, \quext{for some }\,Z_1>2.
\end{equation}
with $c$ independent of $\epsi$. Then, interpolating
between the last two conditions in \eqref{st:n3}, it is clear
that such an exponent can be found if and only if $p>6/5$, as required 
in the statement of Theorem~\ref{teo:main2}.

Next, we notice that, from \eqref{hp:b} and \eqref{st:15}, there follows
\begin{equation}\label{st:a7}
   \| b(\sigma_n,\fhi_n) \|_{L^\infty(0,T;L^p(\Omega))} \le c.
\end{equation}
Let us now write \eqref{nutr} for the approximate sequence and test it by 
$\eta\in W^{1,\infty}(\Omega)$ (as already noted, this condition may be
weakened). Then, using \eqref{st:d12}, \eqref{st:fhi:3d2},
\eqref{st:fhi:3d3} and \eqref{st:a7}, one deduces 
\begin{equation}\label{st:b9b}
  \| \sigma_{n,t} \|_{L^{Z}(0,T;W^{1,R}(\Omega)')} \le c_\epsi,
   \quext{where }\,Z>1,
\end{equation}
and $c_\epsi$ can be taken {\sl independent of }\/ $\epsi$
when either $p>6/5$ or $\sigma_0\in L^q(\Omega)$ for $q>6/5$.

Combining \eqref{st:b9b} with \eqref{st:15} and \eqref{st:fhi:3d3} allows us 
to apply the Aubin-Lions lemma, which gives
\begin{equation}\label{co:n1}
  \sigma_n \to \sigma 
   \quext{strongly in }\,L^2(0,T;W^{1-\iota,S}(\Omega)),
\end{equation}
where $S$ is as in \eqref{st:fhi:3d3} and $\iota>0$ is arbitrarily small.
Here and below, all convergence relations are to be intended to hold
up to the extraction of non-relabelled subsequences of $n\nearrow \infty$.
Note that \eqref{co:n1} implies in particular the pointwise (a.e.) convergence 
$\sigma_n\to \sigma$.

Proving strong convergence of $\fhi_n$ is, of course, much simpler. Using 
\eqref{st:a3}, \eqref{st:11} and \eqref{st:n1}, the Aubin-Lions lemma
actually guarantees 
\begin{equation}\label{co:b1n}
  \fhi_n \to \fhi 
   \quext{strongly in }\,C^0([0,T];H^{1-\iota}(\Omega)) \cap 
    L^{\frac{12p}{12-5p}}(0,T;H^{2-\iota}(\Omega)) \quext{for every }\,\iota>0.
\end{equation}
Again, the above implies in particular almost everywhere
convergence both of $\fhi_n\to \fhi$ and of 
$\nabla\fhi_n\to \nabla\fhi$. This allows us to manage 
the source terms in \eqref{CH1} and 
\eqref{nutr}. Actually, using assumptions \eqref{hp:h}
and \eqref{hp:b} with the first \eqref{st:n3} (notice 
that $5p/3>20/11$ if $p>12/11$), we deduce at least
\begin{align}\label{co:b4}
  & h(\sigma_n,\fhi_n) \to h(\sigma,\fhi)
   \quext{strongly in }\,L^P(Q)~~\text{for every }\,P\in [1,\infty),\\
 \label{co:b5}
  & b(\sigma_n,\fhi_n) \to b(\sigma,\fhi)
   \quext{strongly in }\,L^{20/11}(Q).
\end{align}
We now move to the Korteweg terms in the Brinkman relation \eqref{brink}. 
First of all, notice that, from \eqref{co:b1n}, one deduces in particular
\begin{equation}\label{co:b8}
  \nabla\fhi_n \to \nabla\fhi
   \quext{strongly in }\,L^{4}(0,T;H).
\end{equation}
Combining this fact with the weak convergence resulting from \eqref{st:17},
we then infer
\begin{equation}\label{co:n8}
  \mu_n\nabla\fhi_n \to \mu\nabla\fhi
   \quext{weakly in }\,L^{4/3}(0,T;L^{3/2}(\Omega)).
\end{equation}
Next, let us notice that, combining \eqref{st:a0} with the pointwise 
convergence $\fhi_n\to\fhi$, there follows
\begin{equation}\label{co:n8b}
  \fhi_n \to \fhi
   \quext{strongly in }\,L^{P}(Q)~~%
   \text{for every }\,P\in[1,\infty).
\end{equation}
The above relation and the weak convergence resulting from 
\eqref{st:fhi:3d3} imply
\begin{equation}\label{co:n8n}
  \fhi_n\nabla\sigma_n \to \fhi\nabla\sigma
   \quext{weakly in }\,L^2(0,T;L^S(\Omega)).
\end{equation}
Concerning the cross-diffusion term, as a consequence of the pointwise
convergence of $\sigma_n$ and $\nabla\fhi_n$ and of
\eqref{st:fhi:3d2}, we deduce (at least) 
\begin{equation}\label{co:n9}
  \alpha(\sigma_n)\nabla \fhi_n \to \alpha(\sigma)\nabla \fhi
    \quext{strongly in }\,L^1(Q).
\end{equation}
Finally, we consider the transport terms, both of which are more conveniently
integrated by parts. As already observed, the transport term in \eqref{CH1}
can be treated easily: indeed, using \eqref{st:a2}, the weak $L^2$-convergence 
of $\vu_n$, and the pointwise convergence of $\fhi_n$, there follows
\begin{equation}\label{co:n10}
  \vu_n \fhi_n \to \vu\fhi 
    \quext{weakly in }\,\LDH.
\end{equation}
Concerning the transport term in \eqref{nutr}, a similar argument based on
\eqref{st:d12} (notice that, once more, $Z>1$ is essential)
permits us to obtain at least 
\begin{equation}\label{co:n11}
  \vu_n \sigma_n \to \vu\sigma
    \quext{weakly in }\,L^Z(Q)~~%
    \text{for some }\,Z>1.
\end{equation}
As already observed, for $p\le 6/5$, the argument works only as 
far as $\epsi>0$ is fixed.

The last term we need to manage is the logarithmic one occurring in \eqref{CH2}.
As is customary, to deal with it we will use some tool from the theory
of maximal monotone operators. To this aim, we recall that $\beta(r):=f(r)+\lambda r$ 
denotes the ``monotone part'' of $f$. Hence, of course, it is enough to
prove that $\beta(\fhi_n)\to \beta(\fhi)$ in a suitable sense. 
It is worth observing that, at least formally,
this identification can be obtained just using
the information resulting from the energy estimate. Indeed,
recalling \eqref{st:18} and using that $p$ is {\sl strictly}\/
larger than $1$, one deduces 
\begin{equation}\label{co:n11b}
  \beta(\fhi_n) \to \xi
    \quext{weakly in }\,L^p(Q).
\end{equation}
Then, to identify $\xi$ as $\beta(\fhi)$, one may notice that \eqref{co:n8b} 
holds in particular for $P=p'$ (the conjugate exponent to $p$). Hence,
using the strong-weak closedness of the maximal monotone operator 
induced by $\beta$ and acting from $L^{p'}(Q)$ to $L^p(Q)$ (which 
are reflexive Banach spaces in duality), one readily has
\begin{equation}\label{co:b7}
  \beta(\fhi_n) \to \beta(\fhi)
   \quext{weakly in }\,L^p(Q).
\end{equation}
We notice, however, that adapting this part of the argument to the proposed 
approximation scheme will require some amount of additional work (see the next
section for details).

As a consequence of the above procedure, all nonlinear
terms in the weak formulation \eqref{CH1w}-\eqref{brinkw} pass
to the expected limits. 
To conclude, we observe that, from \eqref{st:15},
\eqref{st:b9b} and a generalized version of the Aubin-Lions
lemma, there follows 
\begin{equation}\label{co:c1}
  \sigma_n \to \sigma 
   \quext{strongly in }\,C^0([0,T];X),
\end{equation}
where $X$ is a negative order Sobolev space such that 
\begin{equation}\label{co:c2}
  L^p(\Omega) \subset \subset X \subset W^{1,\infty}(\Omega)',
\end{equation}
the first inclusion being compact and the second being continuous. 
The existence of such a space is clear as one considers that
$W^{1,\infty}(\Omega)\subset\subset L^{p'}(\Omega)$.
Then, using \eqref{co:c1} and the first \eqref{co:b1n} 
it is immediate to recover the initial data in the sense of \eqref{init}.

\medskip

\noindent%
{\bf End of proof of Theorem~\ref{teo:main2}.}~~%
To conclude the proof, we first have to check the required 
regularity properties \eqref{rego:fhi1}-\eqref{rego:H1}. Actually, most
of them are direct consequences of the ``energy'' bounds \eqref{st:11}-\eqref{st:26},
\eqref{st:a1}, \eqref{st:17} and of the ``entropy'' bounds \eqref{st:ent22}, 
\eqref{st:ent21} (or \eqref{st:ent21b}), and (in the case when we assume $q>p$)
\eqref{st:ent22q}, \eqref{st:ent21q} (or \eqref{st:ent21q4}). 
In particular, the latter conditions imply the third of \eqref{rego:fhi1},
with the exponent $P_0$ specified in \eqref{P0}.

Next, concerning the last regularity condition
in \eqref{rego:fhi1} and the last in \eqref{rego:Ffhi1}, we mimick 
the interpolation argument leading to \eqref{st:n3}, where we can 
take $q$ in place of $p$. This gives 
\begin{equation}\label{st:n3q}
  \| \sigma_n \|_{L^{2}(0,T;L^{\frac{3q}{3-q}}(\Omega))}
   \le c.   
\end{equation} 
Applying to \eqref{CH2} the usual elliptic regularity argument, we then
complete \eqref{rego:fhi1} and \eqref{rego:Ffhi1}. Concerning the 
regularity of $\sigma$, \eqref{rego:sigma1b} follows directly from 
\eqref{st:15} and \eqref{st:ent22} if $q=p$ (i.e., we do not need to 
perform the regularity argument in Subsec.~\ref{subsec:boot})
and from \eqref{st:ent22q} if $q>p$. Regarding \eqref{rego:sigma1a}, 
the second condition follows directly from \eqref{st:fhi:3d3}, while 
the first property is a consequence of \eqref{st:a6} and \eqref{st:d12}; 
in particular, this regularity conditions allows the use of 
$\eta\in W^{1,R}(\Omega)$, with $R$ specified by \eqref{P0},
as a test function in \eqref{nutrw}. 

Next, let us move to the Brinkman-to-Darcy limit.
To this aim, let $(\fhi\ee,\mu\ee,\sigma\ee,\vu\ee)$ be a family 
of weak solutions solving the problem at the level $\epsi$
and let us let $\epsi\searrow0$. Then, we observe that most of 
the estimates are uniform with respect to $\epsi$ and, consequently,
the argument to pass to the limit mostly follows the procedure 
used for letting $n\nearrow \infty$. The only relevant difference 
regarding the fact that, from \eqref{st:14} we can only deduce
that (for a subsequence of $\epsi\searrow0$) 
\begin{equation}\label{st:eeu}
  \vu\ee \to \vu \quext{weakly in }\,\LDH.
\end{equation} 
On the other hand, using \eqref{st:d13}, which holds thanks to the 
assumption $p>6/5$, we can now deduce 
\begin{equation}\label{st:eeu2}
  \sigma\ee \to \sigma \quext{strongly in }\,\LDH.
\end{equation} 
Using \eqref{st:eeu}, \eqref{st:eeu2} and \eqref{st:d13} again, we then 
arrive at 
\begin{equation}\label{st:eeu3}
  \sigma\ee\vu\ee \to \sigma\vu \quext{weakly in }\, L^{Z}(Q),
\end{equation} 
for some $Z>1$. This allows, from one side, to take the limit of the 
transport term in \eqref{nutrw} and, from the other side, to achieve 
the first convergence relation in \eqref{coee:23}, which basically
concludes the proof.

It is finally worth observing that condition \eqref{pos:sigma1} follows 
from \eqref{st:26}, which also implies a uniform in time control
of $\ln\sigma$ in the space of signed measures on $\Omega$.

\medskip

\noindent%
{\bf End of proof of Theorem~\ref{teo:main}.}~~%
In the two-dimensional case, the procedure is simpler and we can rely on the 
stronger estimates obtained in Subsection~\ref{subsec:en1}, which
imply the desired regularity conditions \eqref{rego:fhi}-\eqref{rego:H},
as a direct check shows. In particular, we may notice that 
the last \eqref{rego:fhi} and the last \eqref{rego:Ffhi} are direct consequences
of \eqref{st:ent2d}. Moreover, for fixed $\epsi>0$, combining 
\eqref{rego:vu} with the second \eqref{rego:sigma} and using 
Sobolev's embeddings there follows
\begin{equation}\label{st:2d:tran}
  \| \sigma \vu \|_{L^{2}(0,T;L^{R'}(\Omega))}
   \le c_{R',\epsi},
   \quext{for every }\,R' \in[1,2).
\end{equation} 
Using this piece of information and comparing terms in \eqref{nutr},
it is then easy to get the first \eqref{rego:sigma}, with $R\in(2,\infty]$
being the conjugate exponent to $R'$.

Next, for what concerns the limit $\epsi\searrow 0$, we let as before 
$(\fhi\ee,\mu\ee,\sigma\ee,\vu\ee)$ be a family of weak solutions depending
on $\epsi$. Then, we can mostly proceed as in the three-dimensional
case. The only thing we have to notice regards the first \eqref{coee:13}.
Indeed, by interpolation, \eqref{st:ent12} guarantees
\begin{equation}\label{st:2d:trae}
  \| \sigma\ee \|_{L^{4}(Q)}
   \le c\big( \| \sigma\ee \|_{\LIH}
      + \| \sigma\ee \|_{\LDV} \big) 
  \le c.
\end{equation} 
Combining this with \eqref{coee:14} (which is a consequence of the first 
\eqref{st:14}), we then have
\begin{equation}\label{st:2d:trae2}
  \| \sigma\ee \vu\ee \|_{L^{4/3}(Q)}
   \le c,
\end{equation} 
whence follows the first \eqref{coee:13}. Finally, we observe that,
relation \eqref{nutrw} makes sense for every $\eta\in W^{1,R}(\Omega)$,
$R>2$. This is actually a consequence of \eqref{st:2d:tran}
combined with the fact that \eqref{st:2d:trae} and \eqref{rego:fhi}
imply at least
\begin{equation}\label{st:2d:trae3}
  \| \nabla \sigma \|_{\LDH} + \| \alpha(\sigma) \nabla\fhi \|_{\LDH}
   \le c.
\end{equation} 
Hence the cross-diffusion term in \eqref{nutrw} admits $\eta\in W^{1,R}(\Omega)$ 
as a test function (actually it would also admit $\eta\in V$) and it could be 
written splitting it into its two summands.
\beos\label{2D:meglio}
Adapting the proof of Theorem~\ref{teo:main2}, one may easily check
that, for $d=2$, existence of solutions for $d=2$ can be achieved for 
every under the sole condition \eqref{hp:sigma0} with $p>1$.
Moreover, in this case also the limit Brinkman-to-Darcy
(i.e., $\epsi\searrow0$) can be performed for any $\sigma_0$
complying with \eqref{hp:sigma0} and for any $p>1$
(compare with assumption \eqref{su:q} in the 
three-dimensional case). Actually, for $d=2$ in place of the 
first \eqref{st:n3} we get 
\begin{equation}\label{st:n22d}
  \| \sigma\ee^{p/2} \|_{L^{4}(Q)} 
   \le c \big( \| \sigma\ee^{p/2} \|_{\LIH} + \| \sigma\ee^{p/2} \|_{\LDV} \big)
   \le c,
\end{equation}
where we emphasized the dependence on $\epsi$ in the notation.
Hence, being $p$ is {\sl strictly}\/ larger than $1$, the above implies
\begin{equation}\label{co:n22d}
  \sigma\ee \to \sigma \quext{strongly in }\,L^2(Q), 
\end{equation}
so that, using the weak convergence resulting from the first 
\eqref{st:14} (which is independent of $\epsi$), one deduces
(at least)
\begin{equation}\label{co:n32d}
  \sigma\ee\vu\ee \to \sigma\vu \quext{weakly in }\,L^1(Q),
\end{equation}
which is the key property to achieve the Darcy limit. As noted 
in the introduction, we are planning to analyze in detail 
the two-dimensional case in a forthcoming work.
\eddos


\section{Approximation}
\label{sec:appro}

In this part we sketch a regularization of system \eqref{CH1}-\eqref{incompr}
which may be exploited in order to justify the a priori 
estimates - compactness argument used in the existence
proof. In view of the many nonlinearities involved, 
we will not present the approximation scheme in full detail, as that
would be a rather lengthy and technical procedure. Indeed, we point out 
that, both for the Cahn-Hilliard model and for the chemotaxis
models related to the Keller-Segel system, the literature dealing with 
(analytical and numerical) approximations is very vast: with no
claim of exhaustivity, we refer the reader to \cite{CE,GT,SW,TGG} for 
the Cahn-Hilliard model and to \cite{Sai,SSKHT,TW}
for chemotaxis models. Referring to the quoted papers for more details,
here, we will just focus on the ``new'' aspects of the regularization
arising in connection with our specific 
problem. We also point out that some parts of our argument are somehow reminiscent of 
the procedure devised in \cite{RSS}, to which we refer the reader
for additional considerations. 

That said, we indicate by $n\in\NN$ the regularization parameter, 
intended to go to infinity in the limit, and propose a regularization 
of system \eqref{CH1}-\eqref{incompr} depending on $n$, for which existence
may be proved by standard methods (e.g., a Faedo-Galerkin scheme possibly
complemented by a fixed point argument). The approximation will be
designed so to be compatible with the a-priori estimates obtained 
before. To begin, we need to smooth out the singular
(logarithmic) term $f(\fhi)$ so to ensure the applicability
of some local existence theorem. Actually, in a Faedo-Galerkin
scheme one needs that all the elements of the Galerkin base
belong to the domain of the potential, which is not true if $F$ is 
given by \eqref{Flog} (and so it takes finite values only 
on $[-1,1]$). To regularize $f$, we observe that 
it is sufficient to act on its ``monotone part''
$\beta$. Using some tools from the theory of maximal monotone
operators (cf.\ the monographs \cite{Ba,Br} for some background), 
in Lemma~\ref{lem:linfty} below we will explicitly describe a concrete
example of a ``smooth'' approximating family $\beta_n$,
with controlled growth at infinity, suitably converging to $\beta$ as 
$n\searrow\infty$. Once $\beta_n$ is given, we define
$f_n := \beta_n - \lambda \Id$, $\Id$ denoting the identity function.
Then, noting as $F_n$ a suitable primitive of $f_n$ 
(we may assume $F_n$ to be normalized
so that $F_n(0)=F(0)=0$), one also obtains an approximation of the 
potential.

Once $F_n$ is given, we can introduce a regularized version of the 
energy functional \eqref{energy}, namely
\begin{equation}\label{energyn}
  \calE_n(\fhi,\sigma)
   := \frac12 \| \nabla \fhi \|^2 
   + \frac{1}{2n} \| \Delta \fhi \|^2
   + \io F_n(\fhi) 
   + \io \big( \gammaciapo(\sigma) - \chi \sigma\fhi \big).
\end{equation}
Here we have replaced $F$ with $F_n$ and we have added a further 
regularizing term depending on the Laplacian of $\fhi$. The
reason for such a choice stands in the fact that, as $f$ is smoothed
out, the uniform boundedness of $\fhi$, which is extensively  
used in the estimates (cf.~\eqref{st:a0}) may be lost in the 
approximation. Adding the second-order term, we somehow
compensate this loss of coercivity.

Then, the regularization of the energy, in turn, generates a 
regularized version of the system, which may be stated as
\begin{align}\label{CH1n}
  & \fhi_{n,t} + \vu_n\cdot \nabla \fhi_n - \Delta \mu_n = h(\sigma_n,\fhi_n) - \ell \fhi_n,\\
 \label{CH2n}
  & \mu_n = \frac{1}n \Delta^2 \fhi_n - \Delta \fhi_n + f_n(\fhi_n) - \chi \sigma_n,\\
 \label{nutrn}
  & \sigma_{n,t} + \vu_n \cdot \nabla \sigma_n
   - \Delta\sigma_n + \chi \dive( \alpha(\sigma_n) \nabla \fhi_n ) = b (\sigma_n,\fhi_n),\\
 \label{brinkn}
  & - \epsi \dive (D \vu_n) + \vu_n 
   = \nabla \pi_n + \mu_n \nabla \fhi_n - \chi \fhi_n \nabla \sigma_n,\\
 \label{incomprn}
  & \dive \vu_n = 0.
\end{align} 
Observe that, due to the occurrence of the 
fourth order term in \eqref{CH2n}, an additional boundary
condition is needed. We shall assume
\begin{equation}\label{boundn}
  \dn \Delta \fhi_n = 0
   \quext{on }\,\Gamma\times(0,T).
\end{equation}
Actually, \eqref{CH1n}-\eqref{CH2n} is a particular version of the so-called
``sixth-order Cahn-Hilliard system'', which has been extensively studied in the literature
(we quote, among the other contributions \cite{KNR,SP,SW}) and for which the boundary 
condition \eqref{boundn} is a rather standard choice.

Notice also that there is no need of regularizing the functions $h$ and $b$,
since assumptions \eqref{hp:h} and \eqref{hp:b} already guarantee
sufficient smoothness properties as well as a controlled growth for large values
of $\sigma$ and $\fhi$. We also remark that the expression
of equation \eqref{nutrn} has not been modified compared to \eqref{nutr};
in particular, the equivalent of the ``minimum principle'' \eqref{st:26} still
holds in the approximation. 

We now see how the above regularization impacts on the various parts of the 
existence proof presented in the previous sections.

\smallskip

\noindent%
{\bf Energy estimate.~~}%
We first observe that the balance of mass \eqref{medie} as well as the ``weak''
minimum principle \eqref{stamp} still hold in the approximation. Concerning the energy
bound, the argument outlined in Section~\ref{sec:energy} is fully compatible
with the regularization; however, some differences  arise in the outcome of 
the procedure. To explain this point, we first observe that, as we reproduce the 
estimate by working on \eqref{CH1n}-\eqref{incomprn},
we still obtain the analogue of \eqref{en:14c}. 
Applying Gr\"onwall's lemma, and neglecting for 
simplicity the contribution of the dissipative terms on the \lhs, 
this leads to the bound
\begin{equation}\label{boundn:11}
  \Big\| \calE_n(\fhi_n,\sigma_n) + \| \ln\sigma_n \|_{L^1(\Omega)} \Big\|_{L^\infty(0,T)}
   \le c_T \Big| \calE_n(\fhi_{0,n},\sigma_{0,n}) + \| \ln \sigma_{0,n} \|_{L^1(\Omega)} \Big|.
\end{equation}
Here, $(\fhi_{0,n},\sigma_{0,n})$ are suitable regularizations of the initial
data designed in such a way that the \rhs\ above remains bounded uniformly with 
respect to $n$ and $(\fhi_{0,n},\sigma_{0,n})\to (\fhi_{0},\sigma_{0})$ in a
suitable sense. Actually, one may directly take 
$\sigma_{0,n}\equiv \sigma_0$, while smoothing out $\fhi_0$ is a bit more involved
and can be performed by following the lines of the procedure 
devised in \cite{SP} or referring to other papers dealing with the Cahn-Hilliard
system with singular potential. In particular, one may ensure 
that the analogue of \eqref{hp:fhi0} holds uniformly with respect to $n$
and that 
\begin{equation}\label{boundn:12}
  \| \fhi_{0,n} \|_{H^2(\Omega)}
   \le c n^{1/2},
\end{equation}
so that the regularized initial energy on the \rhs\ of \eqref{boundn:11} is
controlled uniformly in $n$.

The effects of the regularization also have some impact on \eqref{st:11}-\eqref{st:a0}.
In particular, as $F_n$ is defined over the whole real line, there is no reason for 
\eqref{st:a0} to hold. Nevertheless, if we choose $\beta_n$ as in Lemma~\ref{lem:linfty}
below, i.e., we ask that 
\begin{equation}\label{coercn}
  \io F_n(\fhi) \ge \kappa \| \fhi \|_{L^{q_0}(\Omega)}^{q_0} - c,
\end{equation}
for some constants $\kappa>0$ and $c\ge 0$ independent of $n$
and for $q_0>p'$ ($p'$ being the conjugate exponent to $p$),
then the energy functional keeps some coercivity. Indeed,
recalling \eqref{gammaciapo}, one has
\begin{align}\no
  & \io F_n(\fhi) + \io \big( \gammaciapo(\sigma) - \chi \sigma\fhi \big)\\
 \no
  & \mbox{}~~~~~ \ge \kappa \| \fhi \|_{L^{q_0}(\Omega)}^{q_0} 
   + \frac1{p(p-1)} \| \sigma \|_{L^p(\Omega)}^p
   - c \Big( 1 
    + \| \sigma \|_{L^1(\Omega)} + \io |\sigma\ln\sigma|
    + \| \sigma \|_{L^{p}(\Omega)} \| \fhi \|_{L^{p'}(\Omega)} \Big)\\
 \label{coercn2}
   & \mbox{}~~~~~ \ge \frac{\kappa}2 \| \fhi \|_{L^{q_0}(\Omega)}^{q_0} 
   + \frac1{2p(p-1)} \| \sigma \|_{L^p(\Omega)}^p
   - c,
\end{align}
the last inequality following by properly applying Young's inequality and using
that $q_0>p'$. The above relation entails that the approximate energy is {\sl uniformly
coercive}\/ with respect to $n$. As a consequence, when performing the approximate
energy estimate, {\sl at least}\/ a uniform control of $\fhi$ 
in $L^\infty(0,T;L^{q_0}(\Omega))$ is achieved. Moreover, due to the additional 
Laplacian in \eqref{energyn}, \eqref{boundn:11} entails the additional
estimate
\begin{equation}\label{boundn:13}
  \| \fhi_n \|_{L^\infty(0,T;H^2(\Omega))} 
   \le c n^{1/2},
\end{equation}
with the constant $c$ independent of $n$. 

\smallskip

\noindent%
{\bf Entropy estimates.}~~%
We start observing that the first, and simpler, version of the argument, 
as presented in Subsec.~\ref{subsec:en1} for the two-dimensional case,
basically requires no significant variation. On the other hand, for the 
argument presented in Subsec.~\ref{subsec:en2}, some discussion is in order. 
Indeed, the procedure uses in an essential way 
the outcome of the Gagliardo-Nirenberg inequality \eqref{gani:ell}, 
which, in turn, subsumes the uniform boundedness
of $\fhi$ stated in \eqref{st:a0} as a consequence of the choice of the
``singular'' potential \eqref{Flog}, which is however lost in the approximation
as $F$ is smoothed out.
Actually, at the regularized level, the best information we have is given 
by the combination of \eqref{coercn2} and \eqref{boundn:13}. In particular,
the latter relation implies, via Sobolev's embeddings, that, at fixed $n$,
$\fhi_n$ is bounded in the uniform norm; however, at least 
in principle, this bound needs not be uniform with respect to $n$. 

Fortunately, in the next lemma we may prove that $\beta_n$
can be constructed so to ensure the boundedness of $\fhi_n$ independently 
of the approximation parameter.
We believe this property to have an independent interest; hence we present 
it as a lemma so to facilitate future reference. We refer again
to the monographs \cite{Ba,Br} for the terminology as well as to \cite[Chap.~3]{At} 
for what concerns the concepts of graph-convergence and G-convergence.
\bele\label{lem:linfty}
 Let $\Omega$ be a smooth, bounded domain of $\RR^d$, $d\in\{2,3\}$, and 
 $\beta:\RR\to 2^{\RR}$ be a maximal monotone graph, with the\/ {\rm domain}
 $D(\beta)$ satisfying $\overline{D(\beta)}=[-1,1]$. Let $\beta$ be normalized 
 in such a way that $\beta(0)\ni 0$. 
 Let also $q_0\in(2,\infty)$ be a given exponent.
 Then, there exists an approximating family 
 of monotone and locally Lipschitz continuous functions $\beta_n:\RR \to \RR$,
 with $\beta_n(0)=0$, such that, defining 
 \begin{equation}\label{ciapon}
   \betaciapo_n(s) := \int_0^s \beta_n(r)\,\dir,
 \end{equation}
 then it follows that
 \begin{equation}\label{betan11}
   K_n (1 + |s|^{q_0}) \ge \betaciapo_n(s) \ge \kappa |s|^{q_0} - c
    \quext{for every }\,n\in\NN,~s\in \RR,
 \end{equation}
 for suitable constants $\kappa>0$ and $c\ge 0$ (independent of $n$)
 and $K_n>0$ (depending on $n$ and diverging as $n\nearrow\infty$).
 Moreover, for every constant $c_0\ge 0$ there exists a 
 constant ${\ov C}>0$ independent of $n$ and such that, if $\{v_n\}$ 
 is a sequence of functions in $H^2_{\bn}(\Omega)$ satisfying
 \begin{equation}\label{betan12}
   \| v_n \|_V 
    + n^{-\frac12} \| v_n \|_{H^2(\Omega)}
    + \io \betaciapo_n(v_n) \le c_0,
 \end{equation}
 then $\{v_n\}$ also satisfies  
 \begin{equation}\label{betan13}
   \| v_n \|_{L^\infty(\Omega)}\le {\ov C}.
 \end{equation}
 Finally, the family $\beta_n$ converges to $\beta$ in the sense of graphs
 (or G-convergence) in $\RR\times \RR$.
\enle
%
%
\begin{proof}
We denote as $\beta_{0,n}$ the Yosida approximation of $\beta$ of index $n^{-1}$.
Moreover, we define the following sequence of functions:
\begin{equation}\label{indi:11}
   j_n(s):=\begin{cases}
             0 &  \text{~~if }\,s\in[-1,1],\\
             q_0 n^{8q_0} (s-1)^{q_0-1} & \text{~~if }\,s>1,\\
             - q_0 n^{8q_0} |s+1|^{q_0-1} & \text{~~if }\,s<-1.
           \end{cases}
\end{equation}
Using that, for every $s>1$ (respectively, for every $s<-1$) $j_n(s)$
is monotone increasing
(respectively, is monotone decreasing) with respect to $n$ and diverges to 
$+\infty$ (respectively, $-\infty$), it is easy to check that the sequence 
$\{j_n\}$ converges, in the sense of graphs, to the graph $\de I_{[-1,1]}$ (the 
subdifferential of the indicator function of $[-1,1]$). To be precise,
this property may be shown by proving that the antiderivatives of $j_n$
vanishing at $0$ monotonically converge to $I_{[-1,1]}$ 
and applying \cite[Thm.~3.20 and Thm.~3.66]{At}. Then, setting
$\beta_n(s):=\beta_{0,n}(s) + j_n(s)$, we can also see that $\beta_n$ 
converges to $\beta$ in the sense of graphs. In addition, defining 
$\betaciapo_n$ as in \eqref{ciapon}, it is clear that properties
\eqref{betan11} hold. Moreover, one has $\beta_n(0) = 0$ for every
$n\in\NN$ thanks to the properties of Yosida approximations.

Let now \eqref{betan12} hold. Then, computing explicitly the primitive of $j_n$, one 
can easily realize that 
\begin{equation}\label{betan21}
  c_0 \ge \io \betaciapo_n(v_n)
   \ge n^{8q_0} \bigg( \io \big(( v_n - 1)_+\big)^{q_0}
    +  \io \big(( v_n + 1)_-\big)^{q_0} \bigg).
\end{equation}
Let now $T\in C^2(\RR;\RR)$ a truncation operator satisfying the following properties:
\begin{gather}\label{indi:21}
  T(s)\equiv 0~~\text{for all }\,s\in[-1,1], \quad
  T(s)=s-2~~\text{for all }\,s\ge 3, \quad
  T(s)=s+2~~\text{for all }\,s\le - 3,\\
 \label{indi:22}
  | T(s) | \le 1~~\text{for all }\,s\in[-3,3],\quad
   0 \le T'(s) \le 1~~\text{for all }\,s\in\RR, \quad
   | T''(s) | \le c~~\text{for all }\,s\in\RR,
\end{gather}
and for some $c>0$. It is not difficult to construct explicitly a function $T$ satisfying 
all the above properties. Noting that
\begin{equation}\label{betan22}
   \Delta T(v_n) = T'(v_n) \Delta v_n + T''(v_n) | \nabla v_n |^2,
\end{equation}
using Sobolev's embedding and elliptic regularity (recall that 
$v_n\in H^2_{\bn}(\Omega)$) one deduces
\begin{equation}\label{betan22b}
  \| \Delta T(v_n) \| 
   \le c \| \Delta v_n \| + c \| \nabla v_n \|_{L^4(\Omega)}^2
   \le c \big( \| v_n \|_V^2 +  \| \Delta v_n \|^2 \big).
\end{equation}
Hence, recalling \eqref{betan12}
and using again $v_n\in H^2_{\bn}(\Omega)$, there follows 
\begin{equation}\label{betan23}
  \| v_n \|_{H^2(\Omega)}
   \le c \big( \| v_n \| + \| \Delta T(v_n) \| \big)
   \le c_1 n,
\end{equation}
where $c_1>0$ depends on $c_0$ and embedding constants, but is independent of $n$.

Next, let us observe that 
\begin{equation}\label{betan21a}
  | T(v_n) | \le (v_n-1)_+ + (v_n+1)_-
   \quext{a.e.\ in }\,\Omega.
\end{equation}
Hence, from \eqref{betan21} we deduce
\begin{equation}\label{betan21b}
  \| T(v_n) \|_{L^{q_0}(\Omega)}^{q_0} 
   \le c \big( \| (v_n-1)_+ \|_{L^{q_0}(\Omega)}^{q_0}
   + \| (v_n+1)_- \|_{L^{q_0}(\Omega)}^{q_0} \big)
   \le c n^{-8q_0}.
\end{equation}
As a consequence,
\begin{equation}\label{betan21c}
  \| T(v_n) \|
   \le c \| T(v_n) \|_{L^{q_0}(\Omega)}
   \le c n^{-8}.
\end{equation}
Thus, by interpolation (we consider the more difficult case $d=3$,
the two-dimensional setting working with better exponents)
\begin{align}\no
  \| T(v_n) \|_{L^\infty(\Omega)}
  & \le c \| T(v_n) \|_{H^{7/4}(\Omega)}
   \le c \| T(v_n) \|_{H^{2}(\Omega)}^{7/8} \| T(v_n) \|^{1/8}\\
 \label{betan21d}  
  & \le c n^{7/8} n^{-1} = c n^{-1/8}.
\end{align}
Recalling \eqref{indi:21}, the thesis \eqref{betan13} follows.
\end{proof}
With the lemma at disposal, we may notice that, from the approximate
energy estimate \eqref{boundn:11}, there follows in particular that 
\begin{equation}\label{st:12n}
  n^{-\frac12} \| \Delta \fhi_n \|_{\LIH} 
  + \| F_n(\fhi_n) \|_{L^\infty(0,T;L^1(\Omega))}
   \le c.
\end{equation}
Hence, applying the lemma, one obtains
\begin{equation}\label{st:a0n}
  \| \fhi_n \|_{L^\infty(0,T;L^\infty(\Omega))}
   \le c, 
\end{equation}
with $c$ independent of $n$, which allows for the use of the uniform 
boundness of $\fhi$ even when the estimates are performed at the regularized
level. The remainder of the entropy argument seems to require no
further variations.

\smallskip

\noindent%
{\bf Passage to the limit.}~~%
Here, the main point is concerned with the logarithmic term.
Actually, in the previous part (cf., in particular, \eqref{co:b7}), 
the identification of the limit was obtained by relying on
the strong-weak closedness of the maximal monotone operator induced by $\beta$ 
(the monotone part of $f$) from $L^{p'}(Q)$ to $L^p(Q)$ (which are 
reflexive Banach spaces in duality). This property may be stated as
\begin{equation}\label{ws:closed}
  \fhi_n \to \fhi~~\text{strongly in }\,L^{p'}(Q)~\,\text{and }\,%
  \beta(\fhi_n) \to \xi~~\text{weakly in }\,L^p(Q)%
  ~~\Longrightarrow~~\xi\equiv \beta(\fhi),
\end{equation}
the last equality holding in $L^p(Q)$, hence almost everywhere. 

However, as the system is approximated by \eqref{CH1n}-\eqref{incomprn},
some differences arise. First of all, $\beta$ is now replaced by a varying
family $\{\beta_n\}$. In itself, this would not be a problem because the fact
that $\beta_n$ converges to $\beta$ in the sense of graphs (in $\RR$) 
implies, in turn, that the maximal monotone operators
induced by $\beta_n$ G-converge (we refer, again, to \cite[Chap.~3]{At}
for further details) to the the maximal monotone operator induced by $\beta$. 
Such a property would actually imply 
\begin{equation}\label{ws:closedn}
  \fhi_n \to \fhi~~\text{strongly in }\,L^{p'}(Q)~\,\text{and }\,%
  \beta_n(\fhi_n) \to \xi~~\text{weakly in }\,L^p(Q)%
  ~~\Longrightarrow~~\xi\equiv \beta(\fhi).
\end{equation}
There is, however, a further issue which prevents the application of 
the above argument, related to the presence of the fourth-order
elliptic term in \eqref{CH2n}. Indeed, it turns out that the operator
\begin{equation}\label{angle:10}
  B_n(v) := \frac1n \Delta^2 v + \beta_n(v),
\end{equation}
complemented with the boundary conditions \eqref{boundn},
{\sl is not}\/ maximal monotone in the duality between 
$L^{p'}$ and $L^p$. 

To explain this fact, we go back to the Hilbert setting
(indeed, in our modified version of the argument we will be able to 
use Hilbert spaces) and neglect, for simplicity, 
the effects of the time variable. What happens is that,
while the operators $v\mapsto n^{-1} \Delta^2 v$ and $v\mapsto \beta_n(v)$ 
(with proper domains) are {\sl separately}\/ 
maximal monotone in $H$, their sum {\sl is not}.
The reason stants in the lack of the so-called ``angle condition''. 
Roughly speaking, this corresponds to the fact that the $L^2$-scalar 
product of $\Delta^2 v$ and $\beta_n(v)$ has no sign properties. 
In other words, an $L^2$-estimate of the form
\begin{equation}\label{angle:11}
  \big\| n^{-1} \Delta^2 v_n + \beta_n(v_n) \big\| \le c,
\end{equation}
due to the lack of the angle condition,  does {\sl not}\/ imply
\begin{equation}\label{angle:12}
  \big\| n^{-1} \Delta^2 v_n \big\|
    + \| \beta_n(v_n) \| \le c.
\end{equation}
For this reason, the argument in \eqref{ws:closedn} cannot be 
reproduced. To bypass this issue, we go back to the (approximate)
entropy bound and we first recall that \eqref{st:ent21}, as $p>12/11$,
gives at least 
\begin{equation}\label{angle:13}
  \big\| \Delta \fhi_n \big\|_{\LDH} \le c.
\end{equation}
On the other hand, even in the worst case scenario for $p$, the second
\eqref{st:n3} guarantees
\begin{equation}\label{angle:14}
  \| \sigma_n \|_{L^{8/5}(0,T;H)} \le c.
\end{equation}
Hence, using \eqref{st:17} and comparing terms in \eqref{CH2n}, we deduce 
the analogue of \eqref{angle:11}, namely
\begin{equation}\label{angle:15}
  \big\| n^{-1} \Delta^2 \fhi_n + \beta_n(\fhi_n) \big\|_{L^{8/5}(0,T;H)} 
   \le c.
\end{equation}
This estimate, however, as observed above, cannot be decoupled.
For this reason, we need to go back again to \eqref{st:n3} and we 
observe that it also implies
\begin{equation}\label{angle:16}
  \| \sigma_n \|_{L^{2}(0,T;V')}
   \le c \| \sigma_n \|_{L^{2}(0,T;L^{12/7}(\Omega))}
   \le c.
\end{equation}
Then, a further comparison of terms in \eqref{CH2n} gives the estimate
\begin{equation}\label{angle:17}
  \big\| n^{-1} \Delta^2 \fhi_n + \beta_n(\fhi_n) \big\|_{L^2(0,T;V')} \le c,
\end{equation}
which carries a weaker information with respect to space variables,
but, on the other hand, can be decoupled. To see this, we recall
\eqref{co:51}, which, due to the presence of the 
elliptic regularizing term, takes now the form 
\begin{equation}\label{co:51nn}
   \| \Delta\fhi_n \|^2 
   + n^{-1} \| \nabla \Delta\fhi_n \|^2 
   \le c \big( 1 + \| \nabla\mu_n \| \big)
    + \chi \io \nabla \fhi_n \cdot \nabla \sigma_n.
\end{equation}
Then, following the lines of the entropy estimate, we may see that,
in place of \eqref{co:51b}, one has
\begin{equation}\label{co:51bnn} 
   \| \Delta\fhi_n \|^{\frac{12p}{12-5p}}
   + n^{-\frac{6p}{12-5p}} \| \nabla \Delta\fhi_n \|^{\frac{12p}{12-5p}}
   \le c + c \| \nabla\mu_n \|^{\frac{6p}{12-5p}}
    + c \| \sigma_n \|_{L^p(\Omega)}^{\frac{p(5p-4)}{12-5p}}
     \big( 1 + \| \nabla \sigma_n^{p/2} \|^2 \big).
\end{equation}
Hence, using again $\frac{12p}{12-5p}\ge 2$ we additionally deduce 
(at least)
\begin{equation}\label{co:51n2}
   \| n^{-1} \Delta\fhi_n \|_{\LDV} \le c, \quext{or, more 
   precisely, }\,
  \| n^{-1} \Delta^2 \fhi_n \|_{\LDVp} \le c.
\end{equation}
Comparing with \eqref{angle:17}, we then obtain the desired
``decoupled'' information
\begin{equation}\label{angle:18}
  \| \beta_n(\fhi_n) \big\|_{L^2(0,T;V')} \le c.
\end{equation}
Actually, this suffices to pass to the limit provided that we use a 
proper duality argument (an alternative, but in fact equivalent, 
approach may rely on the theory of variational inequalities).
Namely, one may first observe that the sequence of monotone functions 
$\{\beta_n\}$ also induce a family of maximal monotone operators 
from $L^2(0,T;V)$ to $L^2(0,T;V')$. Moreover, such a family
converges, in the sense of G-convergence, 
to the maximal monotone operator induced by 
$\beta$ from $L^2(0,T;V)$ to the power set $2^{L^2(0,T;V')}$:
let us denote by $\beta_w$, i.e., ``weak $\beta$'', the operator
obtained in this way. Actually, as first noticed in \cite{Br2},
$\beta_w$, due to the singular character of $F$, is a
multivalued mapping; in particular there may occur concentration
phenomena on the set where $|\fhi|=1$.

Within this perspective, using the second \eqref{co:b1n} and
\eqref{angle:18}, \eqref{ws:closed} is replaced by
\begin{equation}\label{ws:closedw} 
  \fhi_n \to \fhi~~\text{strongly in }\,L^2(0,T;V)~\,\text{and }\,%
  \beta_n(\fhi_n) \to \xi~~\text{weakly in }\,L^2(0,T;V')%
  ~~\Longrightarrow~~\xi\in \beta_w(\fhi),
\end{equation}
where we used the inclusion symbol because, as said, $\beta_w(\fhi)$ 
can contain more than an element. On the other hand, as 
the identification $\xi\in \beta_w(\fhi)$ is achieved, 
{\sl a posteriori}\/ one can prove that, in fact, $\beta_w(\fhi)$ can be 
replaced by $\beta(\fhi)$. More precisely, the limit $\xi$
of $\beta_n(\fhi_n)$ is a function (and not just a functional),
and it coincides almost everywhere with 
$\beta(\fhi)$. To see this, we first notice that, 
comparing terms in the {\sl limit}\/
equation \eqref{CH2} (where the bi-Laplacian no longer occurs),
one deduces $\xi\in L^{8/5}(0,T;H)$ (the exponent $8/5$ is due to
\eqref{angle:14}). Then, applying, e.g.,
\cite[Prop.~2.5]{BCGG} (see also \cite[Subsec.~2.1]{SP} for 
some additional background) one gets the desired 
(pointwise) equality $\xi=\beta(\fhi)$.

That said, the remainder of the argument used for passing to the limit
with $n\nearrow\infty$ can be adapted to the approximation scheme
up to purely tecnical modifications.


\section*{Acknowledgments}

This work has been partially supported by the PRIN MIUR-MUR Grant 2020F3NCPX 
``Mathematics for industry 4.0 (Math4I4)" and by GNAMPA (Gruppo Nazionale per l'Analisi 
Matematica, la Probabilit\`a e le loro Applicazioni) of INdAM (Istituto Nazionale
di Alta Matematica). The author is grateful to Prof. Kentaro Fujie from 
Tohoku University for useful discussion about the topics of this paper.



\vspace{15mm}

\noindent%
{\bf Author's address:}\\[1mm]
Giulio Schimperna\\
Dipartimento di Matematica, Universit\`a degli Studi di Pavia\\
Via Ferrata, 5,~~I-27100 Pavia,~~Italy\\
E-mail:~~{\tt giulio.schimperna@unipv.it}

\end{document}